\documentclass[12pt,english]{bourbaki}

\usepackage{amsmath}
\usepackage{amsxtra}
\usepackage{amscd}
\usepackage{amsthm}
\usepackage{amsfonts}
\usepackage{amssymb}
\usepackage{eucal}
\usepackage{epsfig}

\newcommand{\thmref}[1]{Theorem~\ref{#1}}
\newcommand{\secref}[1]{\S\ref{#1}}
\newcommand{\defref}[1]{Definition~\ref{#1}}

\newcommand{\propref}[1]{Proposition~\ref{#1}}
\newcommand{\corref}[1]{Corollary~\ref{#1}}
\newcommand{\remref}[1]{Remark~\ref{#1}}
\newcommand{\nc}{\newcommand}
\nc{\cali}{\mathcal}
\nc{\on}{\operatorname}
\nc{\Wick}{{\mb :}}
\nc{\delz}{\partial_z}
\nc{\ddz}{\frac{\partial}{\partial z}}
\nc{\ch}{\mbox{ch}}
\nc{\Z}{{\mbb Z}}
\nc{\C}{{\mbb C}}
\nc{\Oo}{{\cali O}}
\nc{\D}{{\cali D}}
\nc{\cond}{|\,}
\nc{\bib}{\bibitem}
\nc{\pone}{\Pro^1}
\nc{\pa}{\partial}
\nc{\F}{{\cali F}}
\nc{\K}{{\cali K}}
\nc{\arr}{\rightarrow}
\nc{\larr}{\longrightarrow}
\nc{\al}{\alpha}
\nc{\ket}{\rangle}
\nc{\bra}{\langle}
\nc{\W}{{\cali W}}
\nc{\gam}{\bar{\gamma}}
\nc{\Q}{\bar{Q}}
\nc{\q}{\widetilde{Q}}
\nc{\la}{\lambda}
\nc{\ep}{\epsilon}
\nc{\su}{\widehat{{\mf s}{\mf l}}_2}
\nc{\sw}{{\mf s}{\mf l}}
\nc{\g}{{\mf g}}
\nc{\h}{{\mf h}}
\nc{\n}{{\mf n}}
\nc{\ab}{\mf{a}}
\nc{\f}{\widehat{{\cali F}}}
\nc{\is}{{\mb i}}
\nc{\V}{\cali{V}}
\nc{\M}{\widetilde{M}}
\nc{\js}{{\mb j}}
\nc{\bi}{\bibitem}
\nc{\He}{{\cali H}}
\nc{\inv}{^{-1}}
\nc{\vac}{|0\rangle}
\nc{\ol}{\overline}
\nc{\wt}{\widetilde}
\nc{\wh}{\widehat}
\nc{\dst}{\displaystyle}

\nc{\delt}{\partial_t}
\nc{\ddt}{\frac{\partial}{\partial t}}
\nc{\delx}{\partial_x}
\nc{\mb}{\mathbf}
\nc{\mf}{\mathfrak}

\nc{\mbb}{\mathbb}
\nc{\Ctt}{\C((t))}
\nc{\Ct}{\C[t,t\inv]}

\nc{\ghat}{\wh{\g}}
\newcommand\Cx{{\mbb C}^\times}
\nc{\un}{\underline}
\nc{\mc}{\mathcal}
\nc{\BB}{{\mc B}}
\nc{\bb}{{\mf b}}
\nc{\kk}{{\mf k}}
\nc{\Y}{{\mc Y}}
\nc{\frob}{\times}
\nc{\sm}{\setminus}
\nc{\bs}{\backslash}
\nc{\Pp}{{\mathbb P}^1}
\nc{\Aa}{{\mc A}}

\nc{\AutO}{\on{Aut}\Oo}
\nc{\AUTO}{\un{\on{Aut}}\Oo}
\nc{\AUTK}{\un{\on{Aut}}\K}
\nc{\Mcal}{{\mc M}}
\nc{\Heout}{\He_{\out}}
\nc{\Hetil}{{\widetilde\He}}
\nc{\wb}{\overline}

\nc{\Res}{\on{Res}}
\nc{\pitil}{\Pi}
\nc{\Ctil}{\wt{C}}
\nc{\auto}{\on{Aut} \Oo}
\nc{\phitil}{\wt{\phi}}
\nc{\gz}{\g_{\vec z}}
\nc{\tensorM}{\bigotimes_{i=1}^N{\mathbb M}_i}
\nc{\tensorW}{\bigotimes_{i=1}^N W_{\nu_i,k}}
\nc{\B}{{\mc B}}
\nc{\out}{\on{out}}

\nc{\E}{{\mc E}}
\nc{\m}{{\mathfrak m}}

\nc{\gx}{\g^0_{\vec x}}

\nc{\hx}{\He^0_{\vec x}}
\nc{\tensorpi}{\pi_{\nu_1,\ldots,\nu_N}^\kappa}
\nc{\CN}{{\mathcal C}_N}
\nc{\Cn}{{\mathcal C}'_N}
\nc{\Phizw}{\Phi_{\vec w}({\vec z})}
\nc{\Pro}{{\mathbb P}}

\nc{\G}{\wh{\g}}
\nc{\De}{\Delta}

\nc{\us}{\underset}

\nc{\Ll}{\mc L}
\nc{\dR}{\on{dR}}

\nc{\ds}{\displaystyle}
\nc{\T}{{\mc T}}

\nc{\Xn}{\overset{\circ}X{}^n} \nc{\Dn}{\overset{\circ}D{}^n}
\nc{\Dxn}{\overset{\circ}D{}^n_x} \nc{\varphitil}{\wt{\varphi}}

\nc{\lf}{{\mf l}}
\nc{\GL}{{}^L G}
\nc{\Vir}{\on{Vir}}

\date{Juin 2000}
\bbkannee{52\`eme ann\'ee, 1999-2000} 
\bbknumero{875}
\title{Vertex algebras and algebraic curves}

\author{EDWARD FRENKEL}
\address{Department of Mathematics\\
University of California\\
Berkeley, CA 94720, USA}
\email{frenkel@math.berkeley.edu}

\begin{document}
\maketitle


\section{Introduction}

Vertex operators appeared in the early days of string theory as local
operators describing propagation of string states. Mathematical
analogues of these operators were discovered in representation theory
of affine Kac-Moody algebras in the works of Lepowsky--Wilson
\cite{LW} and I. Frenkel--Kac \cite{FK}. In order to formalize the
emerging structure and motivated in particular by the
I. Frenkel--Lepowsky--Meurman construction of the Moonshine Module of
the Monster group, Borcherds gave the definition of vertex algebra in
\cite{Bo}. The foundations of the theory were subsequently laid down
in \cite{FLM,FHL}; in particular, it was shown in \cite{FLM} that the
Moonshine Module indeed possessed a vertex algebra structure. In the
meantime, Belavin, Polyakov and Zamolodchikov \cite{BPZ} initiated the
study of two-dimensional conformal field theory (CFT). Vertex algebras
can be seen in retrospect as the mathematical equivalent of the
central objects of CFT called the chiral symmetry algebras. Moreover,
the key property of associativity of vertex algebras is equivalent to
the property of operator product expansion in CFT, which goes back to
the pioneering works of Polyakov and Wilson. Thus, vertex algebras may
be thought of as the mathematical language of two-dimensional
conformal field theory.

Vertex algebras have turned out to be extremely useful in many areas
of mathematics. They are by now ubiquitous in representation theory of
infinite-dimensional Lie algebras. They have also found applications
in such fields as algebraic geometry, theory of finite groups, modular
functions and topology. Recently Beilinson and Drinfeld have
introduced a remarkable geometric version of vertex algebras which
they called chiral algebras \cite{BD:ch}. Chiral algebras give rise
to some novel concepts and techniques which are likely to have a
profound impact on algebraic geometry.

In this talk we review the theory of vertex algebras with a particular
emphasis on their algebro-geometric interpretation and
applications. We start in \secref{def} with the axiomatic definition
of vertex algebra, which is somewhat different from, but equivalent to
Borcherds' original definition (see \cite{DL,FKRW,Kac vertex}). We
then discuss some of their most important properties and give the
first examples, which come from infinite-dimensional Lie algebras,
such as Heisenberg, affine Kac-Moody and the Virasoro algebra. More
unconventional examples of vertex algebras which are not generated by
a Lie algebra, such as the $\W$--algebras, are reviewed in
\secref{moreex}.

In \secref{coordind} we explain how to make vertex algebras coordinate
independent, thus effectively getting rid of the formal variable
omnipresent in the standard algebraic approach. This is achieved by
attaching to each conformal vertex algebra a vector bundle with a flat
connection on a formal disc, equipped with an intrinsic operation. The
formal variable is restored when we choose a coordinate on the disc;
the independence of the operation from the choice of coordinate
follows from the fact that the group of changes of coordinates is an
``internal symmetry'' of the vertex algebra. Once we obtain a
coordinate independent object, we can give a rigorous definition of
the space of conformal blocks associated to a conformal vertex algebra
and an algebraic curve $X$ (see \secref{confbl}). From the
physics point of view, conformal blocks give rise to ``chiral
correlation functions'' on powers of $X$ with singularities along the
diagonals.

As we vary the curve $X$ and other data on $X$ reflecting the internal
symmetry of a vertex algebra (such as $G$--bundles), the corresponding
spaces of coinvariants, which are the duals of the spaces of conformal
blocks, combine into a sheaf on the relevant moduli space. This sheaf
carries the structure of a (twisted) $\D$--module, as explained in
\secref{sheaves}. One can gain new insights into the structure of
moduli spaces from the study of these $\D$--modules. For instance, one
obtains a description of the formal deformation spaces of the complex
structure or a $G$--bundle on a curve in terms of certain sheaves on
the symmetric powers of the curve \cite{BG,BD,Gi}. Thus, vertex
algebras appear as the local objects controlling deformations of
curves with various extra structures. This raises the possibility that
more exotic vertex algebras, such as $\W$--algebras, also correspond
to some still unknown moduli spaces.

Finally, in \secref{Chiral Algebras} we discuss the relation between
vertex algebras and the Beilinson--Drinfeld chiral algebras. We review
briefly the description of chiral algebras as factorization algebras,
i.e., sheaves on the Ran space of finite subsets of a curve, satisfying
certain compatibilities. Using this description, Beilinson and
Drinfeld have introduced the concept of chiral homology, which can be
thought of as a derived functor of the functor of coinvariants.

The formalism of vertex and chiral algebras appears to be particularly
suitable for the construction of the conjectural geometric Langlands
correspondence between automorphic $\D$--modules on the moduli space
of $G$--bundles on a smooth projective curve $X$ over $\C$ and flat
$^L G$--bundles on $X$, where $G$ is a simple algebraic group and $^L
G$ is the Langlands dual group (see \cite{BD}). The idea is to
construct the automorphic $\D$--modules corresponding to flat $^L
G$--bundles as the sheaves of coinvariants (or, more generally, chiral
homology) of a suitable vertex algebra. We mention below two such
constructions (both due to Beilinson and Drinfeld): in one of them the
relevant vertex algebra is associated to an affine Kac-Moody algebra
of critical level (see \secref{critical}), and in the other it is the
chiral Hecke algebra (see \secref{Chiral Algebras}).

Another application of vertex algebras to algebraic geometry is the
recent construction by Malikov, Schechtman and Vaintrob \cite{MSV} of
a sheaf of vertex superalgebras on an arbitrary smooth algebraic
variety, called the chiral deRham complex, which is reviewed in
\secref{cdo}.

Vertex algebras form a vast and rapidly growing subject, and it is
impossible to cover all major results (or even give a comprehensive
bibliography) in one survey. For example, because of lack of space, I
have not discussed such important topics as the theory of conformal
algebras \cite{Kac vertex,K:talk} and their chiral counterpart,
Lie$^*$ algebras \cite{BD}; quantum deformations of vertex algebras
\cite{Bo3,EK,FR}; and the connection between vertex algebras and
integrable systems.

Most of the material presented below (note that \S\S
\ref{coordind}--\ref{sheaves} contain previously unpublished results)
is developed in the forthcoming book \cite{book}.

I thank A. Beilinson for answering my questions about chiral algebras
and D. Ben-Zvi for helpful comments on the draft of this paper. The
support from the Packard Foundation and the NSF is gratefully
acknowledged.

\section{Definition and first properties of vertex algebras}    \label{def}

\noindent 2.1. Let $R$ be a $\C$--algebra. An $R$--valued {\em formal
power series} (or formal distribution) in variables
$z_1,z_2,\dots,z_n$ is an arbitrary infinite series
$$A(z_1,\dots,z_n)=\sum_{i_1\in\Z} \cdots \sum_{i_n\in\Z}
a_{i_1,\ldots,i_n} z_1^{i_1}\cdots z_n^{i_n},$$ where each
$a_{i_1,\ldots,i_n} \in R$. These series form a vector space denoted
by $R[[z_1^{\pm 1},\dots,z_n^{\pm 1}]]$.  If $P(z_1,\ldots,z_n) \in
R[[z_1^{\pm 1},\ldots,z_n^{\pm 1}]]$ and $Q(w_1,\ldots,w_m) \in
R[[w_1^{\pm 1},\ldots,w_m^{\pm 1}]]$, then their product is a
well-defined element of $R[[z_1^{\pm 1},\ldots,z_n^{\pm 1},w_1^{\pm
1},\ldots,w_m^{\pm 1}]]$. In general, a product of two formal power
series in the same variables does not make sense, but the product of a
formal power series by a polynomial (i.e., a series, such that
$a_{i_1,\ldots,i_n}=0$ for all but finitely many indices) is always
well-defined.

Let $V$ be a $\Z_+$--graded vector space $V=\oplus_{n=0}^\infty V_n$
with finite-dimensional homogeneous components. An endomorphism $T$ of
$V$ is called homogeneous of degree $n$ if $T(V_m) \subset
V_{n+m}$. Denote by $\on{End} V$ the vector space of linear
endomorphisms of $V$, which are finite linear combinations of
homogeneous endomorphisms. This is a $\Z$--graded algebra.

A {\em field} of conformal dimension $\Delta\in\Z_+$ is an $\on{End}
V$--valued formal power series in $z$,
$$\phi(z)=\sum_{j\in\Z}{\phi_j z^{-j-\Delta}}$$ where each $\phi_j$ is
a homogeneous linear endomorphism of $V$ of degree $-j$.  Two fields
$\phi(z)$ and $\psi(z)$ are called mutually {\em local} if there
exists $N \in \Z_+$, such that
\begin{equation}    \label{locality}
(z-w)^N [\phi(z),\psi(w)] = 0
\end{equation}
(as an element of $\on{End} V[[z^{\pm 1},w^{\pm 1}]]$).


Now we can formulate the axioms of vertex algebra.

\begin{defi}\label{VA}

A {\em vertex algebra} is a collection of data:

\begin{enumerate}
\item[$\bullet$] {\em space of states}, a $\Z_+$--graded vector
space $V=\oplus_{n=0}^\infty{V_n}$, with $\on{dim}(V_n)<\infty$;
\item[$\bullet$] {\em vacuum vector} $\vac\in V_0$;
\item[$\bullet$] {\em shift operator} $T:V\to V$ of
degree one;
\item[$\bullet$] {\em vertex operation}, a linear map $Y(\cdot,z):V\to
\on{End} V[[z,z^{-1}]]$ taking each $A\in V_n$ to a field of conformal
dimension $n$.
\end{enumerate}

These data are subject to the following axioms:
\begin{enumerate}
\item[$\bullet$] {\em (vacuum axiom)}
$Y(\vac,z)=\on{Id}_V$. Furthermore, for any $A\in V$ we have
\newline $Y(A,z)\vac\in A + zV[[z]]$ (i.e., $Y(A,z)\vac$ has a
well--defined value at $z=0$, which is equal to $A$).
\item[$\bullet$] {\em (translation axiom)} For any $A\in V$, 
$[T,Y(A,z)]=\partial_z Y(A,z)$ and $T\vac=0$.
\item[$\bullet$] {\em (locality axiom)} All fields $Y(A,z)$ are
mutually local with each other.
\end{enumerate}
\end{defi}

\begin{rema}
The locality axiom was first introduced in \cite{DL}, where it was
shown that it can be used as a replacement for the original axioms of
\cite{Bo,FLM} (see also \cite{Li}).

It is easy to adopt the above definition to the supercase (see, e.g.,
\cite{Kac vertex}). Then $V$ is a superspace, and the above structures
and axioms should be modified
appropriately; in particular, we need to replace the commutator by the
supercommutator in the definition of locality. Then we obtain the
definition of vertex superalgebra. We mostly consider below purely
even vertex superalgebras, but general vertex superalgebras are very
important; for instance, $N=2$ superconformal vertex superalgebras
appear in physical models relevant to mirror symmetry.

The above conditions on $V$ can be relaxed: it suffices to require
that for any $A \in V, v \in V$, we have: $A_n \cdot v = 0$ for $n$
large enough. It is not necessary to require that $\dim V_n < \infty$
and even that $V$ is graded (in that case however one needs to be
careful when dealing with dual spaces). We impose the above stronger
conditions in order to simplify the exposition.

It is straightforward to define homomorphisms between vertex algebras,
vertex subalgebras, ideals and quotients. If $V_1$ and $V_2$ are two
vertex algebras, then $V_1 \otimes V_2$ carries a natural vertex
algebra structure.
\end{rema}

\setcounter{subsection}{1}

\subsection{Example: commutative vertex algebras.}    \label{commex}

Let $V$ be a $\Z$--graded commutative algebra (with finite-dimensional
homogeneous components) with a unit and a derivation $T$ of degree
$1$. Then $V$ carries a canonical structure of vertex algebra. Namely,
we take the unit of $V$ as the vacuum vector $\vac$, and define the
operation $Y$ by the formula
$$
Y(A,z) = \sum_{n\geq 0} \frac{1}{n!} m(T^n A) z^n,
$$
where for $B \in V$, $m(B)$ denotes the operator of
multiplication by $B$ on $V$. It is straightforward to check that all
axioms of vertex algebra are satisfied; in fact, the locality axiom is
satisfied in a strong sense: for any $A, B \in V$, we have:
$[Y(A,z),Y(B,w)] = 0$ (so we have $N=0$ in formula \eqref{locality}).

Conversely, let $V$ be a vertex algebra, in which locality holds in
the strong sense (we call such a vertex algebra {\em commutative}).
Then locality and vacuum axioms imply that $Y(A,z) \in \on{End}
V[[z]]$ for all $A \in V$ (i.e., there are no terms with negative
powers of $z$).
Denote by $Y_A$ the endomorphism of $V$, which is the constant term of
$Y(A,z)$, and define a bilinear operation $\circ$ on $V$ by setting $A
\circ B = Y_A \cdot B.$ By construction, $Y_A Y_B = Y_B Y_A$ for all
$A,B \in V$.  This implies commutativity and associativity of $\circ$
(see, e.g., \cite{Li}). We obtain a commutative and associative
product on $V$, which respects the $\Z$--gradation. Furthermore, the
vacuum vector $\vac$ is a unit, and the operator $T$ is a derivation
with respect to this product. Thus, we see that the notion of
commutative vertex algebra is equivalent to that of $\Z$--graded
commutative associative algebra with a unit and a derivation of degree
$1$.


\begin{rema}
The operator $T$ may be viewed as the generator of infinitesimal
translations on the formal additive group with coordinate
$z$. Therefore a commutative vertex algebra is the same as a
commutative algebra equipped with an action of the formal additive
group. Thus one may think of general vertex algebras as meromorphic
generalizations of commutative algebras with an action of the formal
additive group. This point of view has been developed by Borcherds
\cite{Bo3}, who showed that vertex algebras are ``singular commutative
rings'' in a certain category. He has also considered generalizations
of vertex algebras, replacing the formal additive group by other
(formal) groups or Hopf algebras.

\end{rema}

\subsection{Non-commutative example: Heisenberg vertex algebra.}
\label{ncex}

Consider the space $\C((t))$ of Laurent series in one variable as a
commutative Lie algebra. We define the Heisenberg Lie algebra $\He$ as
follows. As a vector space, it is the direct sum of the space of
formal Laurent power series $\Ctt$ and a one-dimensional space $\C
{\mb 1}$, with the commutation relations
\begin{equation}\label{He cocycle}
[f(t),g(t)]=-\on{Res} f dg {\mb 1}, \quad \quad [{\mb 1},f(t)] = 0. 
\end{equation}
Here $\on{Res}$ denotes the $(-1)$st Fourier coefficient of a Laurent
series. Thus, $\He$ is a one-dimensional central extension of the
commutative Lie algebra $\Ctt$. Note that the relations \eqref{He
cocycle} are independent of the choice of local coordinate $t$. Thus
we may define a Heisenberg Lie algebra canonically as a central
extension of the space of functions on a punctured disc without a
specific choice of formal coordinate $t$.

The Heisenberg Lie algebra $\He$ is topologically generated by the
generators $b_n=t^n, n\in \Z$, and the central element $\mb 1$, and
the relations between them read
\begin{equation}\label{He relations}
[b_n,b_m]=n\delta_{n,-m}{\mb 1}, \quad \quad [{\mb 1},b_n] = 0.
\end{equation}
The subspace $\C[[t]] \oplus {\mb 1}$ is a commutative Lie subalgebra
of $\He$. Let $\pi$ be the $\He$--module induced from the
one-dimensional representation of $\C[[t]] \oplus {\mb 1}$, on which
$\C[[t]]$ acts by $0$ and ${\mb 1}$ acts by $1$. Equivalently, we may
describe $\pi$ as the polynomial algebra
$\C[b_{-1},b_{-2},\dots]$ with $b_n, n<0$, acting by
multiplication, $b_0$ acting by $0$, and $b_n, n>0$, acting as
$\displaystyle{n\frac{\partial}{\partial b_{-n}}}$. The operators
$b_n$ with $n<0$ are known in this context as creation operators,
since they ``create the state $b_n$ from the vacuum $1$'', while the
operators $b_n$ with $n\geq 0$ are the annihilation operators,
repeated application of which will kill any vector in $\pi$. The
module $\pi$ is called the Fock representation of $\He$.

We wish to endow $\pi$ with a structure of vertex algebra. This
involves the following data (\defref{VA}):

\begin{enumerate}
\item[$\bullet$] $\Z_+$ grading: $\deg b_{j_1}\dots b_{j_k} = -\sum_i
{j_i}$.
\item[$\bullet$] Vacuum vector: $\vac=1$.
\item[$\bullet$] The shift operator $T$ defined by the rules: $T\cdot
1=0$ and $[T,b_i]=-ib_{i-1}$. 
\end{enumerate}

We now need to define the fields $Y(A,z)$. To the vacuum vector
$\vac=1$, we are required to assign $Y(\vac,z)=\on{Id}$. The key
definition is that of the field $b(z)=Y(b_{-1},z)$, which generates
$\pi$ in an appropriate sense. Since $\on{deg}(b_{-1})=1$, $b(z)$
needs to have conformal dimension one. We set
$$b(z)=\sum_{n\in\Z}{b_n z^{-n-1}}.$$
The Reconstruction Theorem stated below implies that once we have
defined $Y(b_{-1},z)$, there is (at most) a unique way to extend this
definition to other vectors of $\pi$. This is not very surprising,
since as a commutative algebra with derivation $T$, $\pi$ is freely
generated by $b_{-1}$. Explicitly, the fields corresponding to other
elements of $\pi$ are constructed by the formula
\begin{equation}\label{He fields}
Y(b_{-n_1}b_{-n_2}\dots b_{-n_k},z) = \frac{1}
{(n_1-1)!\dots(n_k-1)!}  \Wick
\delz^{n_1-1}b(z)\dots\delz^{n_k-1}b(z) \Wick.
\end{equation}
The columns in the right hand side of the formula stand for the {\em
normally ordered product}, which is defined as follows. First, let $\;
:\hspace*{-1mm} b_{n_1} \ldots b_{n_k}\hspace*{-1mm}: \;$ be the
monomial obtained from $b_{n_1} \ldots b_{n_k}$ by moving all
$b_{n_i}$ with $n_i<0$ to the left of all $b_{n_j}$ with $n_j \geq 0$
(in other words, moving all ``creation operators'' $b_n, n<0$, to the
left and all ``annihilation operators'' $b_n, n \leq 0$, to the
right). The important fact that makes this definition correct is that
the operators $b_n$ with $n<0$ (resp., $n\geq 0$) commute with each
other, hence it does not matter how we order the creation (resp.,
annihilation) operators among themselves (this property does not hold
for more general vertex algebras, and in those cases one needs to use
a different definition of normally ordered product, which is given
below). Now define $\; :\hspace*{-1mm} \pa_z^{m_1} b(z) \ldots
\pa_z^{m_k} b(z)
\hspace*{-1mm} : \;$ as the power series in $z$ obtained from the
ordinary product $\pa_z^{m_1} b(z) \ldots \pa_z^{m_k} b(z)$ by
replacing each term $b_{n_1} \ldots b_{n_k}$ with $\; :\hspace*{-1mm}
b_{n_1} \ldots b_{n_k} \hspace*{-1mm}: \;$.

\begin{prop}    \label{piva}
The Fock representation $\pi$ with the structure given above satisfies
the axioms of vertex algebra.
\end{prop}

The most difficult axiom to check is locality. The Reconstruction
Theorem stated below allows us to reduce it to verifying locality
property for the ``generating'' field $b(z)$. This is straightforward
from the commutation relations \eqref{He relations}:
\begin{align*}
[b(z),b(w)] &= \sum_{n,m \in \Z} [b_n,b_m] z^{-n-1} w^{-m-1} = \sum_{n
\in \Z} [b_n,b_{-n}] z^{-n-1} w^{n-1} \\ &= \sum_{n \in \Z} n z^{-n-1}
w^{n-1} = \pa_w \delta(z-w),
\end{align*}
where we use the notation $\ds \delta(z-w) = \sum_{m \in \Z} w^m
z^{-m-1}.$

Now $\delta(z-w)$ has the property that $(z-w) \delta(z-w) = 0$. In
fact, the annihilator of the operator of multiplication by $(z-w)$ in
$R[[z^{\pm 1},w^{\pm 1}]]$ equals $R[[w^{\pm 1}]] \cdot \delta(z-w)$
(note that the product of $\delta(z-w)$ with any formal power series
in $z$ or in $w$ is well-defined). More generally, the annihilator of
the operator of multiplication by $(z-w)^N$ in $R[[z^{\pm 1},w^{\pm
1}]]$ equals $\oplus_{n=0}^{N-1} R[[w^{\pm 1}]] \cdot \pa_w^n
\delta(z-w)$ (see \cite{Kac}). This implies in particular that
$(z-w)^2 [b(z),b(w)] = 0$, and hence the field $b(z)$ is local with
itself. \propref{piva} now follows from the Reconstruction Theorem
below.


\subsection{Reconstruction Theorem.}

We state a general result, which provides a
``generators--and--relations'' approach to the construction of vertex
algebras. Let $V$ be a $\Z_+$--graded vector space, $\vac\in V_0$ a
non-zero vector, and $T$ a degree $1$ endomorphism of $V$. Let $S$ be
a countable ordered set and $\{a^\alpha\}_{\al \in S}$ be a collection
of homogeneous vectors in $V$, with $a^{\alpha}$ of degree
$\Delta_\alpha$. Suppose we are also given fields
$$a^\alpha(z)=\sum_{n\in\Z} {a_n^\alpha z^{-n-\Delta_\alpha}}$$ on
$V$, such that the following hold:
\begin{enumerate}
\item[$\bullet$] For all $\alpha$, $a^\alpha(z)\vac \in a^\al + z
V[[z]]$;
\item[$\bullet$] $T\vac=0$ and $[T,a^\alpha(z)]=\delz a^\alpha(z)$ for
all $\alpha$;
\item[$\bullet$] All fields $a^\alpha(z)$ are mutually local;
\item[$\bullet$] $V$ is spanned by lexicographically ordered monomials
$$a^{\alpha_1}_{-\Delta_{\alpha_1}-j_1} \dots
a^{\alpha_m}_{-\Delta_{\alpha_m}-j_m} \vac$$ (including $\vac$), where
$j_1\geq j_2 \geq \ldots \geq j_m \geq 0$, and if $j_i=j_{i+1}$, then
$\al_i\geq \al_{i+1}$ with respect to the order on the set $S$.
\end{enumerate}

\begin{theo}[\cite{FKRW}]    \label{reconstruction}
Under the above assumptions, the above data together with the
assignment
\begin{equation}\label{field formula}
Y(a^{\alpha_1}_{-\Delta_{\alpha_1}-j_1} \dots
a^{\alpha_m}_{-\Delta_{\alpha_m}-j_m}\vac, z)= \frac{1}{j_1!\dots
j_m!} \Wick \delz^{j_1} a^{\alpha_1}(z) \dots \delz^{j_m}
a^{\alpha_m}(z) \Wick,
\end{equation}
define a vertex algebra structure on $V$.
\end{theo}

Here we use the following general definition of the normally ordered
product of fields. Let $\phi(z),\psi(w)$ be two fields of respective
conformal dimensions $\Delta_\phi$, $\Delta_\psi$ and Fourier
coefficients $\phi_n,\psi_n$. The normally ordered product of
$\phi(z)$ and $\psi(z)$ is by definition the formal power series
$$\Wick \phi(z)\psi(z)\Wick=\sum_{n\in\Z}{\left(
\sum_{m\leq-\Delta_\phi}
{\phi_m\psi_{n-m}}+\sum_{m>-\Delta_\phi}{\psi_{n-m}\phi_m}\right)
z^{-n-\Delta_\phi-\Delta_\psi}}.$$ This is a field of conformal
dimension $\Delta_\phi+\Delta_\psi$. The normal ordering of more than
two fields is defined recursively from right to left, so that by
definition $\Wick A(z)B(z)C(z) \Wick=\Wick A(z) (\Wick B(z)
C(z)\Wick)\Wick$. It is easy to see that in the case of the Heisenberg
vertex algebra $\pi$ this definition of normal ordering coincides with
the one given in \secref{ncex}.

\subsection{The meaning of locality.}


The product $\phi(z)\psi(w)$ of two fields is a well--defined
$\on{End} V$--valued formal power series in $z^{\pm 1}$ and $w^{\pm
1}$. Given $v \in V$ and $\varphi \in V^*$, consider the matrix
coefficient
$$
\bra \varphi,\phi(z)\psi(w) v\ket \in \C[[z^{\pm 1},w^{\pm 1}]].
$$
Since $V_n = 0$ for $n < 0$, we find by degree considerations that it
belongs to $\C((z))((w))$, the space of formal Laurent series in $w$,
whose coefficients are formal Laurent series in $z$. Likewise, $\bra
\varphi,\psi(w)\phi(z) v\ket$ belongs to $\C((w))((z))$.

As we have seen in \secref{commex}, the condition that the fields
$\phi(z)$ and $\psi(w)$ literally commute is too strong, and it
essentially keeps us in the realm of commutative algebra. However,
there is a natural way to relax this condition, which leads to the
more general notion of locality. Let $\C((z,w))$ be the field of
fractions of the ring $\C[[z,w]]$; its elements may be viewed as
meromorphic functions in two formal variables. We have natural
embeddings
\begin{equation}\label{two domains}
\C((z))((w))\longleftarrow\C((z,w)) \longrightarrow \C((w))((z)).
\end{equation}
which are simply the inclusions of $\C((z,w))$ into its completions in
two different topologies, corresponding to the $z$ and $w$ axes. For
example, the images of $1/(z-w) \in \C((z,w))$ in $\C((z))((w))$ and
in $\C((w))((z))$ are equal to
$$
\delta(z-w)_- = \sum_{n\geq
0} w^n z^{-n-1}, \quad \quad -\delta(z-w)_+ = -
\frac{1}{z}\sum_{n<0} w^n z^{-n-1},
$$
respectively.

Now we can relax the condition of commutativity of two fields by
requiring that for any $v \in V$ and $\varphi \in V^*$, the matrix
elements $\bra \varphi,\phi(z)\psi(w) v\ket$ and $\bra \varphi,
\psi(w)\phi(z) v\ket$ are the images of the same element
$f_{v,\varphi}$ of $\C((z,w))$ in $\C((z))((w))$ and $\C((w))((z))$,
respectively. In the case of vertex algebras, we additionally require
that $$f_{v,\varphi} \in \C[[z,w]][z^{-1},w^{-1},(z-w)\inv]$$ for all
$v \in V, \varphi \in V^*$, and that the order of pole of
$f_{v,\varphi}$ is universally bounded, i.e., there exists $N \in
\Z_+$, such that $(z-w)^N f_{v,\varphi} \in\C[[z,w]][z^{-1},w^{-1}]$
for all $v,\varphi$. The last condition is equivalent to the condition
of locality given by formula \eqref{locality}.


{} From the analytic point of view, locality of the fields $\phi(z)$
and $\phi(w)$ means the following. Suppose that $\varphi$ is an
element of the restricted dual space $\oplus_{n\geq 0} V_n^*$. Then
$\bra \varphi,\phi(z) \psi(w) v\ket$ converges in the domain
$|z|>|w|$ whereas $\bra \varphi,\psi(w)\phi(z) v\ket$ converges in the
domain $|w|>|z|$, and both can be analytically continued to the same
meromorphic function $f_{v,\varphi}(z,w) \in
\C[z,w][z^{-1},w^{-1},(z-w)^{-1}]$. Then their commutator (considered
as a distribution) is the difference between boundary values of
meromorphic functions, and hence is a delta-like distribution
supported on the diagonal. In the simplest case, this amounts to the
formula $\delta(z-w) = \delta(z-w)_- + \delta(z-w)_+$, which is a
version of the Sokhotsky--Plemelj formula well-known in complex
analysis.

\subsection{Associativity.}    \label{associativ}

Now we state the ``associativity'' property of vertex
algebras. Consider the $V$--valued formal power series
$$
Y(Y(A,z-w)B,w) C = \sum_{n\in\Z} Y(A_n \cdot B,w) C (z-w)^{-n-\Delta_A}
$$
in $w$ and $z-w$. By degree reasons, this is an element of
$V((w))((z-w))$. By locality, $Y(A,z) Y(B,w) C$ is the expansion in
$V((z))((w))$ of an element of $V[[z,w]][z^{-1},w^{-1},(z-w)^{-1}]$,
which we can map to $V((w))((z-w))$ sending $z^{-1}$ to
$(w+(z-w))^{-1}$ considered as a power series in positive powers of
$(z-w)/w$.

\begin{prop}[\cite{FHL,Kac vertex}]    \label{assocprop}
Any vertex algebra $V$ satisfies the following associativity property:
for any $A, B, C \in V$ we have the equality in $V((w))((z-w))$
\begin{equation}    \label{eq:assoc}
Y(A,z) Y(B,w) C = Y(Y(A,z-w)B,w) C.
\end{equation}
\end{prop}

One can also show that for large enough $N \in \Z_+$, $z^N
Y(Y(A,z-w)B,w) C$ belongs to $V[[z,w]][w^{-1},(z-w)^{-1}]$, and its
image in $V[[z^{\pm 1},w^{\pm 1}]]$, obtained by expanding
$(z-w)^{-1}$ in positive powers of $w/z$, equals $z^N Y(A,z) Y(B,w)
C$.

Formula \eqref{eq:assoc} is called the operator product expansion
(OPE). From the physics point of view, it manifests the important
property in quantum field theory that the product of two fields at
nearby points can be expanded in terms of other fields and the small
parameter $z-w$. From the analytic point of view, this formula
expresses the fact that the matrix elements of the left and right hand
sides of the formula, well-defined in the appropriate domains, can be
analytically continued to the same rational functions in $z,w$, with
poles only at $z=0, w=0$, and $z=w$. To simplify formulas we will
often drop the vector $C$ in formula \eqref{eq:assoc}.

For example, in the case of Heisenberg algebra, we obtain:
\begin{equation}    \label{heisOPE}
b(z) b(w) = \frac{1}{(z-w)^2} + \sum_{n\geq 0} \frac{1}{n!} :\pa_w^n
b(w) b(w): (z-w)^n = \frac{1}{(z-w)^2} + :b(z) b(w): .
\end{equation}

Using the Cauchy formula, one can easily extract the commutation
relations between the Fourier coefficients of the fields $Y(A,z)$ and
$Y(B,w)$ from the singular (at $z=w$) terms of their OPE. The result
is
\begin{equation}    \label{commutator}
[A_m,B_k] = \sum_{n>-\Delta_A} \left( \begin{array}{c} m+\Delta_A-1 \\
                        n+\Delta_A-1 \end{array} \right) (A_n \cdot
                        B)_{m+k}.
\end{equation}

\subsection{Correlation functions.}    \label{many}

Formula \eqref{eq:assoc} and the locality property have the following
``multi--point'' generalization. Let $V^*$ be the space of all linear
functionals on $V$.

\begin{prop}[\cite{FHL}]    \label{boot}
Let $A_1,\ldots,A_n \in V$. For any $v \in V, \varphi \in V^*$, and
any permutation $\sigma$ on $n$ elements, the formal power series in
$z_1,\ldots,z_n$,
\begin{equation}    \label{corrfun}
\varphi \left( Y(A_{\sigma(1)},z_{\sigma(1)}) \ldots
Y(A_{\sigma(n)},z_{\sigma(n)}) \vac \right)
\end{equation}
is the expansion in $\C((z_{\sigma(1)})) \ldots ((z_{\sigma(n)}))$ of
an element $f^\sigma_{A_1,\ldots,A_n}(z_1,\ldots,z_n)$ of \newline
$\C[[z_1,\ldots,z_n]][(z_i-z_j)^{-1}]_{i\neq j}$, which satisfies the
following properties: it does not depend on $\sigma$ (so we suppress it
in the notation);
$$
f_{A_1,\ldots,A_n}(z_1,\ldots,z_n) =
f_{Y(A_i,z_i-z_j)A_j,A_1,\ldots,\wh{A_i},\ldots,\wh{A_j},\ldots,
A_n}(z_j,z_1,\ldots,\wh{z_i},\ldots,\wh{z_j},\ldots,z_n)
$$
for all $i\neq j$; and $\pa_{z_i} f_{A_1,\ldots,A_n}(z_1,\ldots,z_n) =
f_{A_1,\ldots,TA_i,\ldots,A_n}(z_1,\ldots,A_n)$.
\end{prop}

Thus, to each $\varphi \in V^*$ we can attach a collection of matrix
elements \eqref{corrfun}, the ``$n$--point functions'' on $\on{Spec}
\C[[z_1,\ldots,z_n]][(z_i-z_j)^{-1}]_{i\neq j}$. They satisfy a
symmetry condition, a ``bootstrap'' condition (which describes the
behavior of the $n$--point function near the diagonals in terms of
$(n-1)$--point functions) and a ``horizontality'' condition. We obtain
a linear map from $V^*$ to the vector space ${\mc F}_n$ of all
collections $\{ f_{A_1,\ldots,A_m}(z_1,\ldots,z_m), A_i \in
V\}_{m=1}^n$ satisfying the above conditions. This map is actually an
isomorphism for each $n\geq 1$. The inverse map ${\mc F}_n \to V^*$
takes $\{ f_{A_1,\ldots,A_m}(z_1,\ldots,z_m), A_i \in V \}_{m=1}^n \in
{\mc F}_n$ to the functional $\varphi$ on $V$, defined by the formula
$\varphi(A) = f_A(0)$. Thus, we obtain a ``functional realization'' of
$V^*$. In the case when $V$ is generated by fields such that the
singular terms in their OPEs are linear combinations of the same
fields and their derivatives, we can simplify this functional
realization by considering only the $n$--point functions of the
generating fields.

For example, in the case of the Heisenberg vertex algebra we consider
for each $\varphi \in \pi^*$ the $n$--point functions
\begin{equation}    \label{npt}
\omega_n(z_1,\ldots,z_n) = \varphi(b(z_1) \ldots b(z_n) \vac).
\end{equation}
By \propref{boot} and the OPE \eqref{heisOPE}, these functions are
symmetric and satisfy the bootstrap condition
\begin{equation}    \label{coordform}
\omega_n(z_1,\ldots,z_n)=\left(
\frac{\omega_{r-2}(z_1,\ldots,\wh{z_i},\ldots,\wh{z_j},\ldots,z_n)}
{(z_i-z_j)^2}+ \mbox{regular} \right).
\end{equation}
Let $\Omega_\infty$ be the vector space of infinite collections
$(\omega_n)_{n\geq 0}$, where $$\omega_n \in
\C[[z_1,\ldots,z_n]][(z_i-z_j)^{-1}]_{i\neq j}$$ satisfy the above
conditions. Using the functions \eqref{npt}, we obtain a
map $\pi^* \to \Omega_\infty$. One can show that this map is an
isomorphism.

\section{More examples}    \label{moreex}

\subsection{Affine Kac-Moody algebras.}


Let $\g$ be a simple Lie algebra over $\C$.  Consider the formal loop
algebra $L\g=\g((t))$ with the obvious commutator. The affine algebra
$\wh{\g}$ is defined as a central extension of $L\g$. As a vector
space, $\wh{\g} = L\g\oplus\C K$, and the commutation relations
read: $[K,\cdot]=0$ and
\begin{equation}\label{Affine relations}
[A\otimes f(t),B\otimes g(t)] =[A,B]\otimes f(t)g(t)
+(\on{Res}_{t=0}fdg(A,B))K
\end{equation}
Here $(\cdot,\cdot)$ is an invariant bilinear form on $\g$ (such a
form is unique up to a scalar multiple; we normalize it as in
\cite{Kac} by the requirement that $(\al_{\on{max}},\al_{\on{max}}) =
2$).

Consider the Lie subalgebra $\g[[t]]$ of $L\g$. If $f,g\in\C[[t]]$,
then $\on{Res}_{t=0}fdg=0$. Hence $\g[[t]]$ is a Lie subalgebra of
$\wh{\g}$. Let $\C_k$ be the one--dimensional representation of
$\g[[t]]\oplus\C K$, on which $\g[[t]]$ acts by $0$ and $K$ acts by
the scalar $k\in\C$. We define the vacuum representation of $\wh{\g}$
of level $k$ as the representation induced from $\C_k$: $\ds
V_k(\g)=\on{Ind}^{\ghat}_{\g[[t]] \oplus \C K} \C_k$.

Let $\{J^a\}_{a=1,\ldots \dim \g}$ be a basis of $\g$. Denote
$J^a_n=J^a\otimes t^n\in L\g$. Then $J^a_n, n\in\Z$, and $K$ form
a (topological) basis for $\wh{\g}$. By the
Poincar\'{e}--Birkhoff--Witt theorem, $V_k(\g)$ has a basis of
lexicographically ordered monomials of the form $J^{a_1}_{n_1}\dots
J^{a_m}_{n_m}v_k$, where $v_k$ is the image of $1 \in \C_k$ in
$V_k(\g)$, and all $n_i <0$. We are now in the situation of the
Reconstruction Theorem, and hence we obtain a vertex algebra structure
on $V_k(\g)$, such that
$$Y(J^a_{-1}v_k,z)=J^a(z):=\sum_{n\in\Z}{J^a_n z^{-n-1}}.$$ The fields
corresponding to other monomials are obtained by formula \eqref{field
formula}. Explicit computation shows that the fields $J^a(z)$ (and
hence all other fields) are mutually local.

\subsection{Virasoro algebra.}    \label{viralg}

The Virasoro algebra $Vir$ is a central extension of the Lie algebra
$\on{Der} \Ctt = \Ctt\ddt$. It has (topological) basis elements
$L_n=-t^{n+1}, n \in \Z$, and $C$ satisfying the relations that
$C$ is central and
$$[L_n,L_m]=(n-m)L_{n+m}+\frac{n^3-n}{12}\delta_{n,-m}C$$ (a
coordinate independent version is given in \secref{ex3}). These
relations imply that the Lie algebra $\on{Der} \C[[t]] = \C[[t]] \ddt$
generated by $L_n, n\geq 0$, is a Lie subalgebra of $Vir$. Consider
the induced representation $\Vir_c=\on{Ind}_{\on{Der}\C[[t]]\oplus\C
C}^{Vir}\C_c$, where $\C_c$ is a one--dimensional
representation, on which $C$ acts as $c$ and $\on{Der} \C[[t]]$ acts
by zero. By the Poincar\'{e}--Birkhoff--Witt theorem, $\Vir_c$ has a
basis consisting of monomials of the form $L_{j_1}\dots L_{j_m}v_c,
j_1\leq j_2\leq j_m\leq -2$, where $v_c$ is the image of $1\in\C_c$ in
the induced representation. By Reconstruction Theorem, we obtain a
vertex algebra structure on $\Vir_c$, such that
$$Y(L_{-2}v_c,z) = T(z) := \sum_{n \in \Z} L_n z^{-n-2}.$$ Note that
the translation operator is equal to $L_{-1}$.

The Lie algebra $\on{Der} \Oo$ generates infinitesimal changes of
coordinates. As we will see in \secref{coordind}, it is important to
have this Lie algebra (and even better, the whole Virasoro algebra)
act on a given vertex algebra by ``internal symmetries''. This
property is formalized by the following definition.

\begin{defi}    \label{confdef}
A vertex algebra $V$ is called {\em conformal}, of central charge
$c\in\C$, if $V$ contains a non-zero {\em conformal vector} $\omega\in
V_2$, such that the corresponding vertex operator
$Y(\omega,z)=\sum_{n\in\Z}{L_n z^{-n-2}}$ satisfies: $L_{-1}=T$,
$L_0|_{V_n}=n \on{Id}$, and $L_2 \omega=\frac{1}{2}c\vac$.
\end{defi}

These conditions imply that the operators $L_n, n \in \Z$, satisfy the
commutation relations of the Virasoro algebra with central charge
$c$. We also obtain a non-trivial homomorphism $\Vir_c \to V$, sending
$L_{-2}v_c$ to $\omega$.

For example, the vector $\frac{1}{2} b_{-1}^2 +\la b_{-2}$ is a
conformal vector in $\pi$ for any $\la \in \C$. The corresponding
central charge equals $1-12\la^2$. The Kac-Moody vertex algebra
$V_k(\g)$ is conformal if $k \neq -h^\vee$, where $h^\vee$ is the dual
Coxeter number of $\g$. The conformal vector is given by the Sugawara
formula $\frac{1}{2(k+h^\vee)} \sum_a (J^a_{-1})^2 v_k$, where $\{ J^a
\}$ is an orthonormal basis of $\g$. Thus, each $\ghat$--module from
the category $\Oo$ of level $k \neq -h^\vee$ is automatically a module
over the Virasoro algebra.

\subsection{Boson--fermion correspondence.}    \label{bfc}

Let ${\mc C}$ be the Clifford algebra associated to the vector space
$\C((t)) \oplus \C((t)) dt$ equipped with the inner product induced by
the residue pairing. It has topological generators $\psi_{n} = t^n,
\psi^*_{n} = t^{n-1} dt, n \in \Z$, satisfying the anti-commutation
relations
\begin{equation}    \label{ferm}
[\psi_{n},\psi_{m}]_+ = [\psi^*_{n},\psi^*_{m}]_+ = 0, \quad \quad
[\psi_{n},\psi^*_{m}]_+ = \delta_{n,-m}.
\end{equation}
Denote by $\bigwedge$ the fermionic Fock representation of ${\mc C}$,
generated by vector $\vac$, such that $\psi_{n} \vac = 0, n\geq 0,
\psi^*_{n} \vac = 0, n>0$. This is a vertex superalgebra with
$$Y(\psi_{-1} \vac,z) = \psi(z) := \ds \sum_{n\in\Z} \psi_{n}
z^{-n-1}, \quad \quad Y(\psi^*_{0} \vac,z) = \psi^*(z) := \ds
\sum_{n\in\Z} \psi_{n}^* z^{-n}.$$

Boson--fermion correspondence establishes an isomorphism between
$\bigwedge$ and a vertex superalgebra built out of the Heisenberg
algebra $\He$ from \secref{ncex}. For $\la \in \C$, let $\pi_\la$ be
the $\He$--module generated by a vector ${\mb 1}_\la$, such that $b_n
\cdot {\mb 1}_\la = \la \delta_{n,0} {\mb 1}_\la, n\geq 0$. For $N \in
\Z_{>0}$, set $V_{\sqrt{N}\Z} = \oplus_{m\in\Z} \pi_{m\sqrt{N}}$. This
is a vertex algebra (resp., superalgebra) for any even $N$ (resp., odd
$N$), which contains $\pi_0 = \pi$ as a vertex subalgebra. The
gradation is given by the formulas $\deg b_n = -n, \deg {\mb 1}_\la =
\la^2/2$ (note that it can be half-integral). The translation operator
is $T = \frac{1}{2} \sum_{n\in\Z} b_n b_{-1-n}$, so $T \cdot {\mb
1}_\la = \la b_{-1} {\mb 1}_\la$. In order to define the vertex
operation on $V_{\sqrt{N}\Z}$, it suffices to define the fields
$Y({\mb 1}_{m\sqrt{N}},z)$. They are determined by the identity
$$\pa_z Y({\mb 1}_\la,z) = Y(T \cdot {\mb 1}_\la,z) = \la :b(z) Y({\mb
1}_\la,z): \quad .$$ Explicitly,
$$
Y({\mb 1}_\la,z) = S_\la z^{\la b_0} \exp \left( - \la \sum_{n<0}
\frac{b_n}{n} z^{-n} \right) \exp \left( - \la \sum_{n>0}
\frac{b_n}{n} z^{-n} \right),
$$
where $S_\la: \pi_\mu \to \pi_{\mu+\la}$ is the shift operator, $S_\la
\cdot {\mb 1}_\mu = {\mb 1}_{\mu+\la}, [S_\la,b_n] = 0, n\neq
0$.

The {\em boson--fermion correspondence} is an isomorphism of vertex
superalgebras $\bigwedge \simeq V_{\Z}$, which maps $\psi(z)$ to
$Y({\mb 1}_{-1},z)$ and $\psi^*(z)$ to $Y({\mb 1}_1,z)$. For more
details, see, e.g., \cite{Kac vertex}.

\subsection{Rational vertex algebras.}    \label{rva}

Rational vertex algebras constitute an important class of vertex
algebras, which are particularly relevant to conformal field theory
\cite{dFMS}. In order to define them, we first need to give the
definition of a module over a vertex algebra.

Let $V$ be a vertex algebra. A vector space $M$ is called a $V$--{\em
module} if it is equipped with the following data:

\begin{enumerate}

\item[$\bullet$] gradation: $M = \oplus_{n \in \Z_+} M_n$;


\item[$\bullet$] operation $Y_M: V \arr \on{End} M[[z,z^{-1}]]$, which
assigns to each $A \in V_n$ a field $Y_M(A,z)$ of conformal dimension
$n$ on $M$;

\end{enumerate}
subject to the following conditions:
\begin{enumerate}

\item[$\bullet$] $Y_M(\vac,z) = \on{Id}_M$;


\item[$\bullet$] For all $A,B \in V$ and $C \in M$, the series
$Y_M(A,z) Y_M(B,w) C$ is the expansion of an element of
$M[[z,w]][z^{-1},w^{-1},(z-w)^{-1}]$ in $M((z))((w))$, and $Y_M(A,z)
Y_M(B,w) C$ \linebreak $= Y_M(Y(A,z-w) B,w) C$, in the sense of
\propref{assocprop}.

\end{enumerate}


A conformal vertex algebra $V$ (see \defref{confdef}) is called {\em
rational} if every $V$--module is completely reducible. It is shown in
\cite{DLM} that this condition implies that $V$ has {\em finitely
many} inequivalent simple modules, and the graded components of each
simple $V$--module $M$ are finite-dimensional. Furthermore, the
gradation operator on $M$ coincides with the Virasoro operator $L^M_0$
up to a shift. Hence we can attach to a simple $V$--module $M$ its
{\em character} $\on{ch} M = \on{Tr}_M q^{L^M_0-c/24}$, where $c$ is
the central charge of $V$. Set $q=e^{2\pi i \tau}$. Zhu \cite{Z1} has
shown that if $V$ satisfies a certain finiteness condition, then the
characters give rise to holomorphic functions in $\tau$ on the
upper-half plane. Moreover, he proved the following remarkable

\begin{theo}    \label{zhu}
Let $V$ be a rational vertex algebra satisfying Zhu's finiteness
condition, and $\{ M_1,\ldots,M_n \}$ be the set of all inequivalent
simple $V$--modules. Then the vector space spanned by $\on{ch} M_i,
i=1,\ldots,n$, is invariant under the action of $SL_2({\mbb Z})$.
\end{theo}

This result has the following heuristic explanation. To each vertex
algebra we can attach the sheaves of coinvariants on the moduli spaces
of curves (see \secref{sheaves}). It is expected that the sheaves
associated to a rational vertex algebra satisfying Zhu's condition are
vector bundles with projectively flat connection. The characters of
simple $V$--modules should form the basis of the space of horizontal
sections of the corresponding bundle on the moduli of elliptic curves
near $\tau=\infty$. The $SL_2({\mbb Z})$ action is just the monodromy
action on these sections.

Here are some examples of rational vertex algebras.

(1) Let $L$ be an even positive definite lattice in a real vector
space $W$. One can attach to it a vertex algebra $V_L$ in the same way
as in \secref{bfc} for $L=\sqrt{N}\Z$. Namely, we define the
Heisenberg Lie algebra associated to $W_{\C}$ as the Kac-Moody type
central extension of the commutative Lie algebra $W_{\C}((t))$. Its
Fock representation $\pi_W$ is a vertex algebra. The vertex structure
on $\pi_W$ can be extended to $V_L = \pi_W \otimes \C[L]$, where
$\C[L]$ is the group algebra of $L$ \cite{Bo,FLM}. Then $V_L$ is
rational, and its inequivalent simple modules are parameterized by
$L'/L$ where $L'$ is the dual lattice \cite{D1}. The corresponding
characters are theta-functions. The vertex algebra $V_L$ is the chiral
symmetry algebra of the free bosonic conformal field theory
compactified on the torus $W/L$.

Note that if we take as $L$ an arbitrary integral lattice, then $V_L$
is a vertex superalgebra.

(2) Let $\g$ be a simple Lie algebra and $k \in \Z_+$. Let $L_k(\g)$
be the irreducible quotient of the $\ghat$--module $V_k(\g)$. This is
a rational vertex algebra, whose modules are integrable
representations of $\ghat$ of level $k$ \cite{FZ}. The corresponding
conformal field theory is the Wess-Zumino-Witten model (see
\cite{dFMS}).

(3) Let $L_{c(p,q)}$ be the irreducible quotient of
$\on{Vir}_{c(p,q)}$ (as a module over the Virasoro algebra), where
$c(p,q)=1-6(p-q)^2/pq$, $p,q>1, (p,q)=1$. This is a rational vertex
algebra \cite{Wa}, whose simple modules form the ``minimal model'' of
conformal field theory defined by Belavin, Polyakov, and Zamolodchikov
\cite{BPZ} (see also \cite{dFMS}).

(4) The Moonshine Module vertex algebra $V^{\natural}$ (see below).

Zhu \cite{Z1} has attached to each vertex algebra $V$ an associative
algebra $A(V)$, such that simple $V$--modules are in one-to-one
correspondence with simple $A(V)$--modules.

\subsection{Orbifolds and the Monster.}


Let $G$ be a group of automorphisms of a vertex algebra $V$. Then the
space $V^G$ of $G$--invariants in $V$ is a vertex subalgebra of
$V$. When $G$ is a finite group, {\em orbifold theory} allows one to
construct $V^G$--modules from twisted $V$--modules. For a finite order
automorphism $g$ of $V$, one defines a $g$--twisted $V$--module by
modifying the above axioms of $V$--module in such a way that
$Y_M(A,z)$ has monodromy $\la^{-1}$ around the origin, if $g \cdot A =
\la v$ (see \cite{DLM}). Then conjecturally all $V^G$--modules can be
obtained as $V^G$--submodules of $g$--twisted $V$--modules for $g \in
G$.

Now suppose that $V$ is a {\em holomorphic} vertex algebra, i.e.,
rational with a unique simple module, namely itself. Let $G$ be a
cyclic group of automorphisms of $V$ of order $n$ generated by
$g$. Then the following pattern is expected to hold (see \cite{D2}):
(1) for each $h \in G$ there is a unique simple $h$--twisted
$V$--module $V_h$; (2) $V_h$ has a natural $G$--action, so we can
write $V_h = \oplus_{i=0}^{n-1} V_h(i)$, where $V_h(i) = \{ v \in V_h|
g \cdot v = e^{2\pi i k/n} v \}$; then each $V_h(i)$ is a simple
$V^G$--module, and these are all simple $V^G$--modules; (3) The vector
space $\wt{V} = \oplus_{h \in G} V_h(0)$ carries a canonical vertex
algebra structure, which is holomorphic.

A spectacular example is the construction of the {\em Moonshine
Module} vertex algebra $V^{\natural}$ by I. Frenkel, Lepowsky and
Meurman \cite{FLM}. In this case $V$ is the vertex algebra $V_\Lambda$
associated to the Leech lattice $\Lambda$ (it is self-dual, hence
$V_\Lambda$ is holomorphic), and $g$ is constructed from the
involution $-1$ on $\Lambda$. Then $V^\natural = V_\Lambda(0) \oplus
V_{\Lambda,g}(0)$ is a vertex algebra, whose group of automorphisms is
the Monster group \cite{FLM}. Moreover, $V^\natural$ is holomorphic
\cite{D2}. Conjecturally, $V^\natural$ is the unique holomorphic
vertex algebra $V$ with central charge $24$, such that $V_1 = 0$. Its
character is the modular function $j(\tau) - 744$. More generally, for
each element $x$ of the Monster group, consider the Thompson series
$\on{Tr}_{V^\natural} x q^{L_0-1}$. The Conway--Norton conjecture, in
the formulation of \cite{FLM}, states that these are Hauptmoduls for
genus zero subgroups of $SL_2({\mbb R})$. This was first proved in
\cite{FLM} in the case when $x$ commutes with a certain involution,
and later in the general case (by a different method) by R. Borcherds
\cite{Bo1}.


\subsection{Coset construction.}    \label{coset}

Let $V$ be a vertex algebra, and $W$ a vector subspace of $V$. Denote
by $C(V,W)$ the vector subspace of $V$ spanned by vectors $v$ such
that $Y(A,z) \cdot v \in V[[z]]$ for all $A \in W$.  Then $C(V,W)$ is
a vertex subalgebra of $V$, which is called the {\em coset vertex
algebra} of the pair $(V,W)$. In the case $W=V$, the vertex algebra
$C(V,V)$ is commutative, and is called the {\em center} of $V$. An
example of a coset vertex algebra is provided by the
Goddard--Kent--Olive construction \cite{GKO}, which identifies
$C(L_1(\sw_2) \otimes L_k(\sw_2),L_{k+1}(\sw_2))$ with
$L_{c(k+2,k+3)}$ (see \secref{rva}). For other examples, see
\cite{dFMS}.

\subsection{BRST construction and $\W$--algebras.}    \label{Walg}

For $A \in V$, let $y(A) = \on{Res} Y(A,z)$. Then associativity
implies: $[y(A),Y(B,z)] = Y(y(A) \cdot B,z)$, and so $y(A)$ is an
infinitesimal automorphism of $V$. Now suppose that $V^\bullet$ is a
vertex (super)algebra with an additional $\Z$--gradation and $A \in
V^1$ is such that $y(A)^2=0$. Then $(V^\bullet,y(A))$ is a complex,
and its cohomology is a graded vertex (super)algebra. Important
examples of such complexes are provided by {\em topological vertex
superalgebras} introduced by Lian and Zuckerman \cite{LZ}. In this
case the cohomology is a graded commutative vertex algebra, but it has
an additional structure of Batalin-Vilkovisky algebra.

Another example is the BRST complex of quantum hamiltonian
reduction. We illustrate it in the case of quantum Drinfeld-Sokolov
reduction \cite{FF:gd}, which leads to the definition of
$\W$--algebras. Let $\g$ be a simple Lie algebra of rank $\ell$, and
$\n$ its upper nilpotent subalgebra.

Let ${\mc C}$ be the Clifford algebra associated to the vector space
$\n((t)) \oplus \n^*((t)) dt$ equipped with the inner product induced
by the residue pairing. It has topological generators $\psi_{\al,n} =
e_\al \otimes t^n , \psi^*_{\al,n} = e_\al^* \otimes t^{n-1} dt, \al
\in \De_+, n \in \Z$, satisfying the anti-commutation relations
\eqref{ferm}.
Let $\bigwedge_{\n}$ be its Fock representation, generated by vector
$\vac$, such that $\psi_{\al,n} \vac = 0, n\geq 0, \psi^*_{\al,n} \vac
= 0, n>0$. This is a vertex superalgebra which is the tensor product
of several copies of the vertex superalgebra $\bigwedge$ from
\secref{bfc}.
Introduce an additional $\Z$--gradation on ${\mc C}$ and
$\bigwedge_{\n}$ by setting $\deg \psi^*_{\al,n} = - \deg \psi_{\al,n}
= 1, \deg \vac = 0$.

Now consider the vertex superalgebra $C^\bullet_k(\g) = V_k(\g)
\otimes \bigwedge_{\n}^\bullet$. It carries a canonical differential
$d_{\on{st}}$ of degree $1$, the standard differential of
semi-infinite cohomology of $\n((t))$ with coefficients in $V_k$
\cite{Fe,FGZ}; it is equal to the residue of a field from
$C^\bullet_k(\g)$ (see \cite{FF:gd}). Let $\chi = \sum_{i=1}^\ell
\on{Res} \psi^*_{\al_i}(z)$ be the Drinfeld-Sokolov character
\cite{DS} of $\n((t))$. Then $d = d_{\on{st}} + \chi$ is a
differential on $C^\bullet_k(\g)$, and the cohomology
$H^\bullet_k(\g)$ of this differential is a vertex superalgebra.

\begin{theo} $H^0_k(\g)$ is a vertex algebra generated by elements
$W_i$ of degrees $d_i+1, i=1,\ldots,\ell$, where $d_i$ is the $i$th
exponent of $\g$ (in the sense of the Reconstruction Theorem), and
$H^i_k(\g) = 0, i\neq 0$.
\end{theo}

This theorem is proved in \cite{FF:gd,FF:im} for generic $k$ and in
\cite{dBT} for all $k$. The vertex algebra $H^0_k(\g)$ is called the
$\W$--{\em algebra} associated to $\g$ and denoted by $\W_k(\g)$. We
have: $\W_k(\sw_2) = \on{Vir}_{c(k)}$, where $c(k) = 1 -
6(k+1)^2/(k+2)$. For $\g=\sw_3$, the $\W$--algebra was first
constructed by Zamolodchikov, and for $\g=\sw_N$ by Fateev and
Lukyanov \cite{FL}. Since $V_k(\g)$ is conformal for $k\neq -h^\vee$
(see \secref{viralg}), $\W_k(\g)$ is also conformal, with $W_1$
playing the role of conformal vector. On the other hand,
$\W_{-h^\vee}(\g)$ is a commutative vertex algebra, which is
isomorphic to the center of $V_{-h^\vee}(\g)$ \cite{FF:gd} (see
\secref{critical}). The simple quotient of $\W_k(\g)$ for $k=-h^\vee +
p/q$, where $p,q$ are relatively prime integers greater than or equal
to $h^\vee$, is believed to be a rational vertex algebra. Moreover, if
$\g$ is simply-laced and $q=p+1$, this vertex algebra is conjecturally
isomorphic to the coset algebra of $(L_1(\g) \otimes
L_{p-2}(\g),L_{p-1}(\g))$, see \cite{FKW}.

For any $V_k(\g)$--module $M$, the cohomology of the complex $(M
\otimes \bigwedge_{\n}^\bullet,d)$ is a $\W_k(\g)$--module. This
defines a functor, which was studied in \cite{FKW}.

The $\W$--algebras exhibit a remarkable duality: $\W_k(\g) \simeq
\W_{k'}({}^L\g)$, where $^L\g$ is the Langlands dual Lie algebra to
$\g$ and $(k+h^\vee) r^\vee = (k'+{}^L h^\vee)^{-1}$ (here $r^\vee$
denotes the maximal number of edges connecting two vertices of the
Dynkin diagram of $\g$), see \cite{FF:gd}. In the limit $k \to
-h^\vee$ it becomes the isomorphism of \thmref{center}.

\section{Coordinate independent description of vertex algebra
structure}    \label{coordind}

Up to this point, we have discussed vertex algebras in the language of
formal power series. The vertex operation is a map $Y:V\to \on{End}
V[[z,z\inv]]$, or, equivalently, an element of the completed tensor
product $V^*((z)) \widehat{\otimes} \on{End} V$. Thus, we can view
$Y$ as an $\on{End} V$--valued section of a vector bundle over the
punctured disc $D^\times = \on{Spec} \C((z))$ with fiber $V^*$. The
question is whether one can define this bundle in such a way that this
section is canonical, i.e., independent of the choice of $z$. In order
to do that, we need a precise description of the transformation
properties of $Y$. Once we understand what type of geometric object
$Y$ is, we will be able to give a global geometric meaning to vertex
operators on arbitrary curves.

\subsection{The group $\on{Aut} \Oo$.}
Denote by $\Oo$ the complete topological $\C$--algebra $\C[[z]]$, and
let $\AutO$ be the group of continuous automorphisms of $\Oo$. Such an
automorphism is determined by its action on the generator
$z\in\C[[z]]$. Thus as a set this group consists of elements
$\rho(z)\in\C[[z]]$ of the form $\rho_1 z + \rho_2 z^2 + \ldots$, with
$\rho_1 \in \Cx$, endowed with the composition law
$(\rho*\mu)(z)=\mu(\rho(z))$. It is easy to see that $\AutO$ is a
proalgebraic group, $\underset{\longleftarrow}{\on{lim}} \on{Aut}
\C[[z]]/(z^n)$. It is the semi-direct product ${\mathbb G}_m \ltimes
\on{Aut}_+{\cali O}$, where $\on{Aut}_+{\cali O}$ is the subgroup that
consists of the transformations of the form $z \to z+a_2z^2+\cdots$,
and ${\mathbb G}_m$ is the group of rescalings $z \to az, a\neq
0$. Furthermore, $\on{Aut}_+{\cali O}$ is a prounipotent group. The
Lie algebra of $\on{Aut} \Oo$ is $\on{Der}_0 \Oo = z\C[[z]] \pa_z$,
which is a semi--direct product of $\C\pa_z$ and $\on{Der}_+ \Oo =
z^2\C[[z]] \pa_z$.  Any representation of $\on{Der}_0 \Oo$, on which
$z \pa_z$ is diagonalizable with integral eigenvalues and $\on{Der}_+
\Oo$ acts locally nilpotently, can be exponentiated to a representation
of $\AutO$.

Given a vertex algebra and a field $Y(A,z)$, it makes sense to
consider a new field $Y(A,\rho(z))$. We now seek an action of
$\on{Aut} {\cali O}$ on $V$, $A \mapsto R(\rho) \cdot A$, such that
$Y(A,\rho(z))$ is related to $Y(R(\rho)A,z)$ in some reasonable
way. In the theory of vertex algebras, actions of $\AutO$ usually
arise from the action of the Virasoro algebra. Namely, let $V$ be a
conformal vertex algebra (see \defref{confdef}). Then the Fourier
coefficients $L_n$ of $Y(\omega,z)$ satisfy the commutation relations
of the Virasoro algebra with central charge $c$. The operators $L_n,
n\geq 0$, then define an action of the Lie algebra $\on{Der}_+ \Oo$
($L_n$ corresponds to $-t^{n+1} \pa_t$). It follows from the axioms of
vertex algebra that this action can be exponentiated to an action of
$\AutO$. For $f(z) \in \AutO$, denote by $R(f):V \to V$ the
corresponding operator.

Given a vector field $\ds v = \sum_{n\geq -1}{v_n z^{n+1} \pa_z}$, we
assign to it the operator $\ds {\mathbf v}=-\sum_{n\geq -1}{v_n L_n}$.
{}From formula \eqref{commutator} we obtain:
\begin{equation}\label{infinitesimal action}
[{\mathbf v},Y(A,w)]=-\sum_{m\geq -1}{\frac{1}{(m+1)!}(\pa^{m+1}_w v(w))
Y(L_m A,w)}.
\end{equation}
By ``exponentiating'' this formula, we obtain the following identity,
due to Y.-Z. Huang \cite{Hu}:
\begin{equation}    \label{Hformula}
Y(A,t) = R(\rho)Y(R(\rho_t)\inv A,\rho(t))R(\rho)\inv,
\end{equation}
for any $\rho \in \on{Aut} \C[[z]]$ (here $\rho_t(z) = \rho(t+z) -
\rho(t)$, considered as formal power series in $z$ with coefficients
in $\C[[t]]$).

\subsection{Vertex algebra bundle.}    \label{global construction}

Now we can give a coordinate-free description of the operation $Y$.
Let $X$ be a smooth complex curve. Given a point $x \in X$, denote by
$\Oo_x$ the completion of the local ring of $x$, by ${\mf m}_x$ its
maximal ideal, and by $\K_x$ the field of fractions of $\Oo_x$. A
formal coordinate $t_x$ at $x$ is by definition a topological
generator of ${\mf m}_x$.
Consider the set of pairs $(x,t_x)$, where $x \in X$ and $t_x$ is a
formal coordinate at $x$. This is the set of points of a scheme $\wh
X$ of infinite type, which is a principal $\on{Aut} \Oo$--bundle over
$X$. Its fiber at $x \in X$ is the $\AutO$--torsor $\Ll_x$ of all
formal coordinates at $x$. Given a
finite-dimensional $\on{Aut}\Oo$--module $V$, let ${\cali V}={\wh
X}\underset{\on{Aut}\Oo}{\times} V$ be the vector bundle associated to
$V$ and $\wh{X}$. The fiber of $\V$ at $x\in X$ is the $\Ll_x$--twist
of $V$, $\V_x = \Ll_x \underset{\AutO}\times V$.

A conformal vertex algebra is an $\AutO$--module, which has a
filtration $V_{\leq i} := \oplus_{k=0}^i V_k$ by finite-dimensional
$\AutO$--submodules. We obtain the directed system $(\V_{\leq i})$ of
the corresponding vector bundles of finite rank on $X$ and embeddings
$\V_{\leq i} \hookrightarrow \V_{\leq j}, i\leq j$. We will denote
this system simply by $\V$, thinking of it as an inductive limit of
bundles of finite rank. Likewise, by $\V^*$ we will understand the
inverse system of bundles $(\V_{\leq i})^*$ and surjections $(\V_{\leq
j})^* \twoheadrightarrow (\V^i)^*$, thinking of it as a projective
limit of bundles of finite rank.

The $\AutO$--bundle $\wh{X}$ on $X$ above carries an action of the Lie
algebra $\on{Der} \Oo$, which is compatible with the $\AutO$--action
and simply transitive, i.e., $\on{Der} \Oo \wh{\otimes} \Oo_{\wh{X}}
\simeq \Theta_{\wh{X}}$. Since $V$ is a $\on{Der} \Oo$--module, by the
general construction of \secref{gener} we obtain a flat connection
$\nabla:\V\to\V\otimes \Omega$ on $\V$. Locally, $\nabla$ can be
written as $d+L_{-1}\otimes dz$, where $d$ the deRham differential.
The connection $\nabla$ on $\V$ gives us a connection $\nabla^*$ on
the dual bundle $\V^*$.
To be precise, $\V$ is not an inductive limit of flat bundles, since
$\nabla(\V_{\leq i}) \subset \V_{i+1} \otimes \Omega$, and likewise,
$\V^*$ is not a projective limit of flat bundles (it is a $\D$--module
on $X$, which is not quasicoherent as an $\Oo$--module).

Now let us restrict $\V^*$ to the punctured disc $D_x^\times =
\on{Spec} \K_x$ around $x \in X$. We want to define an $\on{End}
\V_x$--valued meromorphic section ${\mc Y}_x$ of $\V^*$ on
$D^\times_x$. Pick a formal coordinate $z$ at $x$. With this choice,
identify ${\cali V}_x$ with $V$ and trivialize $\V|_{D_x}$. We define
our section ${\mc Y}_x$ in this trivialization through its matrix
coefficients: to each triple $c\in V$, $c'\in V^*$ and a section
$F:D_x \to V$ we need to assign a meromorphic function on the disc,
which is $\C$--linear in $c$ and $c'$ and ${\cali O}$--linear in
$F$. It suffices to assign such a function to triples $c,c',F$, where
$F$ is the constant section equal to $A \in V$ with respect to our
trivialization. We set this function to be equal to $\langle c' |
Y(A,z)|c \rangle$.

\begin{theo}    \label{universal section}
The section $\Y_x$ defined this way is canonical, i.e., independent of
the choice of formal coordinate $z$ at $x$. Moreover, $\Y_x$ is
horizontal: $\nabla^* \Y_x=0$.
\end{theo}

Thus, we see that a conformal vertex algebra gives rise to a vector
bundle on any smooth curve and a canonical horizontal section of the
restriction of this bundle to the neighborhood of each point with
values in the endomorphisms of the corresponding twist of $V$. The
fact that this section is canonical is equivalent to formula
\eqref{Hformula}.
Flatness follows from the formula $\pa_z Y(B,z)=Y(L_{-1}B,z)$, which
holds in any vertex algebra.

\begin{rema}    \label{multpoints}
Let $\V^{\boxtimes n}$ be the $n$--fold external tensor power of
$\V$. Consider the $n$--point functions \eqref{corrfun}. According to
\propref{boot}, these are elements of
$\C[[z_1,\ldots,z_n]][(z_i-z_j)^{-1}]_{i\neq j}$. Trivializing the
restriction of $\V^{\boxtimes n}$ to $D_x^n = \on{Spec}
\C[[z_1,\ldots,z_n]]$ using the coordinate $z$, we obtain a
$\V_x^*$--valued section of $(\V^{\boxtimes n})^*|_{D_x^n}$ with poles
on the diagonals. Then this section is canonical, i.e., independent of
the choice of $z$, and horizontal.
\end{rema}

\subsection{Examples.}    \label{ex2}
If $W$ is an $\AutO$--stable subspace of $V$, it gives rise to a
subbundle ${\mc W}$ of $\V$. The universal section ${\mc Y}_x$ then
gives rise to an $\on{End} \V_x$--valued section of ${\mc
W}^*|_{D_x^\times}$. This way one obtains many familiar geometric
objects.

Suppose $A\in V_\Delta$ satisfies $L_n \cdot A=0, L_0 A=\Delta A$.
Then subspace $\C A\subset V$ is a one-dimensional $\on{Aut} \Oo$
submodule of $V$. The line bundle associated to this module is nothing
but the line bundle $\Omega^{-\Delta}$ of $\Delta$--differentials on
$X$. Thus we obtain a line subbundle ${\cali L}_A$ of $\V$. Our
section ${\mc Y}_x$ therefore gives us an $\on{End} \V_x$--valued
section of ${\cali L}_A^{-1} = \Omega^{\Delta}$. In other words, the
$\on{End} \V_x$--valued $\Delta$--differential $Y(A,z) (dz)^{\Delta}$
does not depend on the choice of formal coordinate $z$ at $x$.

Next, consider the Heisenberg vertex algebra $V=\pi$, with conformal
structure given by the vector $\frac{1}{2} b_{-1}^2 + b_{-2}$. For
each smooth curve $X$, $\pi$ gives rise to a vector bundle that we
denote by $\Pi$. The first piece $\pi_0$ of our filtration on $\pi =
\C[b_n]_{n<0}$ is the one-dimensional subspace spanned by the vacuum
vector $1$. This is a trivial representation of the group $\on{Aut}
\Oo$. Hence it gives rise to a trivial subbundle $\Oo_X\subset
\Pi$. The next piece in the filtration, $\pi_{\leq 1} = \C 1 \oplus \C
b_{-1}$, gives rise to a rank two subbundle of $\Pi$, which we denote
by ${\cali B}$. It dual bundle is an extension
\begin{equation}    \label{Bla}
0 \to \Omega \to {\cali B}^* \to \Oo_X \to 0.
\end{equation}
Our universal section ${\mc Y}_x$ gives rise to an $\on{End}
\Pi_x$--valued section of ${\cali B}^*|_{D_x^\times}$, which projects
onto the section $\on{End} \Pi_x$--valued section of $\Oo_X$ equal to
$\on{Id}$ (this follows from the vacuum axiom). Explicit computation
shows that the space of sections of ${\cali B}^*$ which project onto
the section $1$ of $\Oo_X$ (equivalently, splittings of \eqref{Bla})
is canonically isomorphic to the space of {\em affine connections}
(they may also be described as affine structures, see \cite{Gu}). Thus,
we conclude that the field $b(z)$ (or, more precisely, the expression
$\pa_z + b(z)$) transforms as an $\on{End} \Pi_x$--valued affine
connection on $D_x^\times$.  Similarly, we obtain connections on
$G$--bundles from the subspace $V_k(\g)_{\leq 1}$ of the Kac-Moody
vertex algebra $V_k(\g)$.

Let $V$ be a conformal vertex algebra with central charge $c$. It
follows from the definition that the vector space $\C\vac\oplus
\C\omega$ is $\on{Aut}\Oo$--stable. Hence it gives rise to a rank two
subbundle ${\cali T}_c$ of $\V$. Its dual bundle is an extension
\begin{equation}    \label{Tc}
0 \to \Omega^2 \to{\cali T}_c^*\to \Oo_X \to 0.
\end{equation}
The universal section ${\mc Y}_x$ gives rise to an $\on{End}
\V_x$--valued section of ${\cali T}_c^*$ over the punctured disc,
which projects onto $\on{Id} \in \on{End} \V_x \otimes
\Oo_x$. Explicit calculation shows that the space of sections of
${\cali T}_c^*$ projecting onto $1\in\Gamma(\Oo_X)$ (equivalently,
splittings of \eqref{Tc}) is isomorphic to the space of self--adjoint
second order differential operators $\rho:\Omega^{-\frac{1}{2}}
\to\Omega^{\frac{3}{2}}$ with constant symbol $\dfrac{c}{6}$. Locally,
such an operator can be written as
$\frac{c}{6}\delz^2+q(z)$. Operators of this form with symbol $1$ are
known as {\em projective connections} (they may also be described as
projective structures, see \cite{Gu}). Thus, we see that for $c\neq 0$
the field $T(z)$ (called the stress tensor in conformal field theory),
or more precisely, the expression $\pa_z^2 + \frac{6}{c} T(z)$
transforms as an $\on{End} \V_x$--valued projective connection on the
punctured disc.

\subsection{General twisting property.}    \label{gentw}

We have seen in \secref{global construction} that the vertex operation
$Y$ is in some sense invariant under the $\AutO$--action (see formula
\eqref{Hformula}). Therefore $Y$ gives rise to a well-defined
operation on the twist of $V$ by any $\AutO$--torsor. This ``twisting
property'' follows from the fact that $\AutO$ acts on $V$ by
``internal symmetries'', that is by exponentiation of Fourier
coefficients of the vertex operator $Y(\omega,z)$. It turns out that
vertex algebras exhibit a similar twisting property with respect to
any group ${\mc G}$ of internal symmetry, i.e., a group (or more
generally, an ind-group) obtained by exponentiation of Fourier
coefficients of vertex operators. Using associativity, one can obtain
an analogue of formula \eqref{Hformula} for the transformations of $Y$
under the action of ${\mc G}$. This formula means that we get a
well-defined operation on the twist of $V$ by any ${\mc G}$--torsor.

For example, let $V$ be a vertex algebra with a {\em
$\ghat$--structure} of level $k\neq -h^\vee$, i.e., one equipped with
a homomorphism $V_k(\g) \to V$, whose image is not contained in
$\C\vac$. Recall that $V_k(\g)$ is a conformal vertex algebra if
$k\neq -h^\vee$ (see \secref{viralg}). Given such a homomorphism
$V_k(\g) \to V$, the Fourier coefficients of the corresponding fields
$Y(\omega,z)$ and $Y(J^a,z)$ generate an action of the Lie algebra
$Vir \ltimes \ghat$ on $V$. Suppose that the action of its Lie
subalgebra $\on{Der}_0 \Oo \ltimes \g(\Oo)$ can be exponentiated to an
action of the group $\AutO \ltimes G(\Oo)$ on $V$ (then we say that
this group acts on $V$ by internal symmetries). Let ${\mc P}$ be a
$G$--bundle on a smooth curve $X$. Denote by $\wh{\mc P}$ the
principal $\on{Aut}\Oo \ltimes G(\Oo)$--bundle over $X$, whose fiber
at $x$ consists of pairs $(z,s)$, where $z$ is a formal coordinate at
$x$ and $s$ is a trivialization of ${\mc P}|_{D_x}$. Let $\V^{\mc P}$
be the $\wh{\mc P}$--twist of $V$. Then the vertex operation $Y$ on
$V$ gives rise to a canonical section $\Y_x^{\mc P}$ of $(\V^{\mc
P})^*|_{D_x^\times}$ with values in $\on{End} \V^{\mc P}_x$.


\section{Conformal blocks}    \label{confbl}

\noindent 5.1. The results of the previous section allow us to assign
to a conformal vertex algebra a vector bundle with connection on any
smooth curve $X$, and a family of local structures on it.
We can now associate to any compact curve $X$ an invariant of the
vertex algebra structure.

\begin{defi}    \label{confblocks}
Let $V$ be a conformal vertex algebra, $X$ a smooth projective curve,
and $x$ a point of $X$. A linear functional $\varphi$ on $\V_x$ is called
a {\em conformal block} if $\varphi(\Y_x \cdot A) \in
\Gamma(D_x^\times,\V^*)$ can be extended to a regular section of
$\V^*$ on $X\bs x$ for all $A\in\V_x$.

The set of conformal blocks is a vector subspace of $\V_x^*$, denoted
by $C(X,x,V)$.
\end{defi}

Let $\varphi \in C(X,x,V)$, and $A \in \V_x$. Denote by $\varphi_A$
the corresponding section of $\V^*$ over $X\bs x$. If $W$ is an
$\AutO$--stable subspace of $V$, then $\varphi_A$ can be projected
onto a section of $\W^*$ over $X\bs x$.  Note that according to the
vacuum axiom, $\varphi_{\vac}$ is actually a regular section of
$\V^*$ over the whole $X$, and so is its projection. In particular,
taking as $W$ the two-dimensional subspace of $V$ spanned by $\vac$
and $\omega$, we assign to each conformal block a projective
connection on $X$ (if $c\neq 0$). Denote $\varphi(\cdot \vac)$ as
$\bra \cdot \ket_\varphi$. If we choose a formal coordinate $z$ at
$x$, then we can write this projective connection as $\pa_z^2 +
\frac{6}{c} \bra T(z) \ket_\varphi$. If $c=0$, then we obtain a
quadratic differential $\bra T(z) \ket_\varphi dz^2$.

\begin{rema}    \label{twconf}
Let $V$ be a vertex algebra with $\ghat$--structure of level $k\neq
-h^\vee$, and ${\mc P}$ a $G$--bundle on $X$. A conformal block {\em
twisted by} ${\mc P}$ is by definition a linear functional $\varphi$
on $\V^{\mc P}_x$, such that $\varphi(\Y^{\mc P}_x \cdot A)$ (see
\secref{gentw}) can be extended to a regular section of $(\V^{\mc
P})^*$ on $X\bs x$ for all $A\in\V^{\mc P}_x$. We denote the
corresponding space by $C^{\mc P}(X,x,V)$.
\end{rema}

\setcounter{subsection}{1}

\subsection{Examples.}    \label{ex3}
In the case of the Kac-Moody vertex algebra $V_k(\g)$ the
\defref{confblocks} can be simplified. In this case the subspace $(\g
\otimes t^{-1}) v_k$ is $\AutO$--invariant, and therefore we have a
surjection $(\V_k(\g))^* \twoheadrightarrow \g^* \otimes
\Omega$. Hence for each $\varphi \in C(X,x,\pi), A \in \V_x$, we
obtain a regular $\g^*$--valued one-form on $X \bs x$, whose
restriction to $D_x^\times$ is $\varphi(J(z) \cdot A) dz$, where $J(z)
= \sum_a J_a \otimes J^a(z)$.

Now observe that any $f = \sum f_n z^n \in \g \otimes \K_x$ gives rise
to an operator on $\V_k(\g)_x$, $\wt{f} = \on{Res}_x (f,J(z)) dz =
\sum f_n b_n$, which does not depend on the choice of the coordinate
$z$. Since by definition $\varphi(J(z) \cdot A) dz$ extends to a
regular one-form on $X \bs x$, we obtain from the residue theorem that
$\varphi(\wt{f} \cdot A) = 0$ for all $f \in \g_{x,\out} = \g \otimes
\C[X\bs x]$. Hence we obtain a map from $C(X,x,V_k(\g))$ to the space
of $\g_{x,\out}$--invariant functionals on $\V_k(\g)_x$.

\begin{lemm} The space of conformal blocks $C(X,x,V_k(\g))$ is
isomorphic to the space of $\g_{x,\on{out}}$--invariant functionals on
$\V_k(\g)_x$.
\end{lemm}

Thus we recover the common definition of conformal blocks as
$\g_{x,\out}$--invariants. They have been extensively studied recently
because of the relation to the moduli spaces of $G$--bundles on curves
(see \secref{globalstr} below).  More generally, we have: $C^{\mc
P}(X,x,V_k(\g)) =$ $\on{Hom}_{\g^{\mc P}_{x,\out}}(\V_k(\g)^{\mc
P}_x,\C)$, where $\g^{\mc P}_{x,\out}$ is the Lie algebra of sections
of ${\mc P} \us{G}\times \g$ over $X\bs x$.

The situation becomes more subtle in the case of the Virasoro vertex
algebra $\Vir_c$. The analogue of $\g_{x,\out}$ is then the Lie
algebra $\on{Vect}(X\bs x)$ of vector fields on $X\bs x$. By analogy
with the Kac-Moody case, we would like to define a homomorphism
$\on{Vect}(X\bs x) \to \on{End}{\mc V}ir_{c,x}$ sending $\xi(z) \pa_z
\in \on{Vect}(X\bs x)$ to $\on{Res}_x \xi(z) T(z) dz$. But it is not
clear whether this map is independent of the choice of formal
coordinate $z$, unless $c=0$ (when the field $T(z)$ transforms as a
quadratic differential). According to \secref{ex2}, if $c \neq 0$,
then $\pa_z^2 + \frac{6}{c}T(z)$ is naturally a section of a rank two
bundle ${\mc T}_c^*$. Hence it can be paired with sections of the rank
two bundle ${\mc T}_c \otimes \Omega$ (we could split this bundle by
choosing a projective structure on $X$, but we prefer not to do
that). One can show that the sheaf $\wh{\mc T}_c = ({\mc T}_c \otimes
\Omega)/d\Oo$, which is an extension of $\Theta$ by $\Omega/d\Oo$, has
a canonical Lie algebra structure (see \cite{Wi,BS}), such that $Vir_x
= \Gamma(D_x^\times,\wh{\mc T}_c)$ is isomorphic to the Virasoro
algebra and acts on ${\mc V}ir_{c,x}$. The map $\Gamma(X\bs x,\wh{\mc
T}_c) \to Vir_x$ factors through $\on{Vect}(X\bs x)$, and hence the
above map $\on{Vect}(X\bs x) \to \on{End}{\mc V}ir_{c,x}$ is
well-defined. One can then show that the space of conformal blocks
$C(X,x,\on{Vir}_c)$ is isomorphic to the space of $\on{Vect}(X\bs
x)$--invariant functionals on ${\mc V}ir_{c,x}$.

\subsection{General invariance condition.}    \label{coinv}

For general vertex algebras the analogue of the above invariance
condition is defined as follows. Let $\V \stackrel{\nabla}{\arr} \V
\otimes \Omega$, where $\V \otimes \Omega$ is placed in degree $0$, be
the de Rham complex of the flat vector bundle $\V$, considered as a
complex of sheaves on the curve $X$. Here $\nabla$ is the connection
operator defined in \secref{global construction} (recall that $\nabla
= d + L_{-1} \otimes dz$). For $\Sigma \subset X$, denote by
$H^0_{\on{dR}}(\Sigma,\V \otimes \Omega)$ the degree $0$ cohomology of
the restriction of this complex to $\Sigma$. If $\Sigma$ is affine or
$D_x^\times$, then $H^0_{\on{dR}}(\Sigma,\V \otimes \Omega)$ is simply
the quotient of $\Gamma(\Sigma,\V \otimes \Omega)$ by the image of
$\nabla$.

There is a linear map $\gamma: \Gamma(D_x^\times,\V \otimes \Omega)
\arr \on{End} \V_x$ sending $\mu \in \Gamma(D_x^\times,\V \otimes
\Omega)$ to a linear operator $O_\mu=\Res_x\bra \Y_x,\mu\ket$ on
$\V_x$. Since the residue of a total derivative vanishes, $\gamma$
factors through $U(\V_x) = H^0_{\on{dR}}(D_x^\times,\V \otimes
\Omega)$. If we choose a formal coordinate $z$ at $x$, then we can
identify $U(\V_x)$ with $U(V)$, the completion of the span of all
Fourier coefficients of all vertex operators $Y(A,z), A \in
V$. According to formula \eqref{commutator}, $U(\V_x)$ is a Lie
algebra. Note that in the case of Kac-Moody and Virasoro vertex
algebras, $U(\V_x)$ contains $\wh{\g}_x$ and $Vir_x$, respectively, as
Lie subalgebras.

Given an open $\Sigma \subset X$, we have a canonical map
$H^0_{\on{dR}}(\Sigma,\V \otimes \Omega) \rightarrow U(\V_x)$. Denote
its image by $U_\Sigma(\V_x)$. It follows from the Beilinson-Drinfeld
chiral algebra formalism (see \secref{Chiral Algebras}) that
$U_\Sigma(\V_x)$ is a Lie subalgebra of $U(\V_x)$.  The image of
$U_\Sigma(\V_x)$ in $\on{End} \V_x$ is the span of all operators
$O_\mu$ for $\mu\in\Gamma(\Sigma,\V\otimes\Omega)$. The residue
theorem implies:

\begin{prop}
The space of conformal blocks $C(X,x,V)$ coincides with the space of
$U_{X\bs x}(\V_x)$--invariant functionals on $\V_x$.
\end{prop}

Therefore the dual space to the space of conformal blocks is the {\em
space of coinvariants} $H(X,x,V) = \V_x/ U_{X\bs x}(\V_x) \cdot \V_x$.

The definition of the spaces $C(X,x,V)$ and $H(X,x,V)$ can be
generalized to the case of multiple (distinct) points $x_1,\dots x_n$
on the curve $X$, at which we ``insert'' arbitrary conformal
$V$--modules $M_1,\ldots,M_n$. A $V$--module $M$ is called conformal
if the gradation operator on $M$ coincides, up to a shift, with the
Fourier coefficient $L_0^M$ of the field $Y_M(\omega,z)$. Assume for
simplicity that the eigenvalues of $L_0^M$ on all modules $M_i$ are
integers. Then each $M_i$ is an $\AutO$--module, and so we can attach
to it a vector bundle ${\mc M}_i$ on $X$ with a flat connection (a
general conformal $V$--module gives rise to a vector bundle with
connection on the space of pairs $(x,\tau_x)$, where $x \in X$ and
$\tau_x$ is a non-zero tangent vector to $X$ at $x$).

The space of conformal blocks $C_V(X,(x_i),(M_i))_{i=1}^n$ is then by
definition the space of linear functionals $\varphi$ on
$\Mcal_{1,x_1}\otimes\cdots\otimes{\Mcal}_{n,x_n}$, such that for any
$A_i\in{\Mcal}_{i,x_i}, i=1,\dots,n$, there exists a section of
${\V}^*$ over $X\setminus\{x_1,\dots,x_n\}$, whose restriction to each
${D}_{x_i}^\times$ equals
$$\varphi(A_1\otimes\cdots\otimes \Y_{M_i}\cdot A_i\otimes \cdots
\otimes A_n).$$ Informally, one can say that the local sections of
$\V^*$ over the discs around the points, obtained by acting with
vertex operators at those points, can be ``glued together'' into a
single meromorphic section of $\V^*$. There is a canonical
isomorphism
\begin{equation}    \label{propagation}
C_V(X,(x_1,\ldots,x_n,y),(M_1,\ldots,M_n,V)) \simeq
C_V(X,(x_1,\ldots,x_n),(M_1,\ldots,M_n)),
\end{equation}
given by $\varphi\mapsto \varphi|_{M_1\otimes\cdots\otimes
M_N\otimes\vac}$.

\begin{rema} Let $V$ be a rational vertex algebra (see \secref{rva}).
Then we obtain a functor from the category of pointed algebraic curves
with insertions of simple $V$--modules at the points to the category
of vector spaces, which sends $(X,(x_i),(M_i))$ to
$C_V(X,(x_i),(M_i))_{i=1}^n$. This is a version of {\em modular
functor} corresponding to $V$ \cite{Se,Gaw,Z2}. It is known that
in the case of rational vertex algebras $L_k(\g)$ or $L_{c(p,q)}$ the
spaces of conformal blocks are finite-dimensional. Moreover, the
functor can be extended to the category of pointed stable curves with
insertions. It then satisfies a factorization property, which
expresses the space of conformal blocks associated to a singular curve
with a double point in terms of conformal blocks associated to its
normalization. The same is expected to be true for general rational
vertex algebras.
\end{rema}



\subsection{Functional realization.}    \label{functreal}

The spaces $C_V(X,(x_i),(M_i))_{i=1}^n$ with varying points
$x_1,\ldots,x_n$ can be organized into a vector bundle on $\Xn :=
X^n\backslash \Delta$ (where $\Delta$ is the union of all diagonals)
with a flat connection, which is a subbundle of $({\mc M}_1 \boxtimes
\ldots \boxtimes {\mc M}_n)^*$.
Consider for example the spaces $C_V(X,(x_i),(V))_{i=1}^n$. By
\eqref{propagation}, all of them are canonically isomorphic to each
other and to $C_V(X,x,V)$. Hence the corresponding bundles of
conformal blocks are canonically trivialized. For $\varphi \in
C_V(X,x,V)$, let $\varphi_n$ denote the corresponding element of
$C_V(X,(x_i),(V))_{i=1}^n$. Let $A_i(x_i)$ be local
sections of $\V$ near $x_i$. Evaluating our conformal blocks on them,
we obtain what physicists call the chiral correlation functions
\begin{equation}    \label{phin}
\varphi_{n}(A_1(x_1) \otimes \ldots \otimes A_n(x_n)) \sim
\bra A_1(x_1) \ldots A_n(x_n) \ket_\varphi.
\end{equation}

Moreover, we can explicitly describe these functions when $x_i$'s are
very close to each other, in the neighborhood of a point $x \in X$: if
we choose a formal coordinate $z$ at $x$ and denote the coordinate of
the point $x_i$ by $z_i$, then
\begin{equation}    \label{product}
\varphi_n(A_1(z_1) \otimes \ldots \otimes A_n(z_n)) =
\varphi(Y(A_1,z_1) \ldots Y(A_n,z_n)\vac).
\end{equation}
We obtain that the restriction of the $n$--point chiral correlation
function to $\Dxn = D_x^n\backslash \Delta :=$ $\on{Spec}
\C[[z_1,\ldots,z_n]][(z_i-z_j)^{-1}]_{i\neq j}$ coincides with the
$n$--point function corresponding to the functional $\varphi$
introduced in \secref{many}.

According to \remref{multpoints}, the matrix elements given by the
right hand side of formula \eqref{product} give rise to a horizontal
section of $(\V^{\boxtimes n})^*$ on $\Dxn$. Therefore by
\propref{boot} we obtain an embedding $\V_x^* \hookrightarrow
\Gamma_\nabla(\Dxn,(\V^{\boxtimes n})^*)$, where $\Gamma_\nabla$
stands for the space of horizontal sections.  If $\varphi \in \V_x^*$
is a conformal block, then \defref{confblocks}
implies that the image of $\varphi \in C_V(X,x,V)$ in
$\Gamma_\nabla(\Dxn,(\V^{\boxtimes n})^*)$ extends to a horizontal
section $\varphi_n$ of $(\V^{\boxtimes n})^*$ over $\Xn$. Thus, we
obtain an embedding $C_V(X,x,V) \hookrightarrow
\Gamma(\Xn,(\V^{\boxtimes n})^*) = \Gamma_\nabla(X^n,j_* j^*
(\V^{\boxtimes n})^*)$, where $j: \Xn \hookrightarrow X^n$.

Imposing the bootstrap conditions on the diagonals from
\propref{boot}, we obtain a quotient $\V_n^*$ of $j_* j^*
(\V^{\boxtimes n})^*$ (the precise definition of $\V_n^*$ uses the
operation ${\mc Y}^{(2)}$ introduced in \thmref{invariant OPE}). Since
by construction our sections satisfy the bootstrap conditions, we have
embeddings $\V_x^* \hookrightarrow \Gamma_\nabla(D_x^n,\V_n^*)$ and
$C_V(X,x,V) \hookrightarrow \Gamma_\nabla(X^n,\V_n^*)$. A remarkable
fact is that these maps are isomorphisms for all $n>1$. Thus, using
correlation functions we obtain {\em functional realizations} of the
(dual space of) a vertex algebra and its spaces of conformal blocks.

\begin{rema} Note that neither $(\V^{\boxtimes n})^*$ nor
$\V_n^*$ is quasicoherent as an $\Oo_{X^n}$--module. It is more
convenient to work with the dual sheaves $\V^{\boxtimes n}$ and
$\V_n$, which are quasicoherent. The sheaf $\V_n$ on $X^n$ is defined
by Beilinson and Drinfeld in their construction of factorization
algebra corresponding to $\V$ (see \secref{Chiral Algebras}). Note
that the space $\Gamma_\nabla(X,\V_n^*)$ is dual to the top deRham
cohomology $H_{\on{dR}}^n(X,\V_n \otimes \Omega_{X^n})$ (which is
therefore isomorphic to $H(X,x,V)$).
\end{rema}


In the case of Heisenberg algebra, the functional realization can be
simplified, because we can use the sections of the exterior products
of the canonical bundle $\Omega$ instead of the bundle $\Pi^*$ (here
we choose $\frac{1}{2} b_{-1}^2$ as the conformal vector, so that
$b(z) dz$ transforms as a one-form). Namely, as in \secref{many} we
assign to $\varphi \in \Pi_x^*$ the collection of polydifferentials
$$
\varphi(b(z_1) \ldots b(z_n)\vac) dz_1 \ldots dz_n.
$$
These are sections of $\Omega^{\boxtimes n}(2\Delta)$ over $D_x^n$,
which are symmetric and satisfy the bootstrap condition
\eqref{coordform}. Let $\Omega_\infty$ be the sheaf on $X$, whose
sections over $U \subset X$ are collections of sections of
$\Omega^{\boxtimes n}(2\Delta)$ over $U^n$, for $n\geq 0$, satisfying
these conditions (the bootstrap condition makes sense globally because
of the identification $\Omega^{\boxtimes 2}(2\Delta)/\Omega^{\boxtimes
2}(\Delta) \simeq \Oo$). Then we have canonical isomorphisms $\Pi^*_x
\simeq \Gamma(D_x,\Omega_\infty), C(X,x,\pi)\simeq
\Gamma(X,\Omega_\infty)$.

The same construction can be applied to the Kac-Moody vertex algebras,
see \cite{BD:poly}.

\begin{rema}
If we consider conformal blocks with general $V$--module insertions,
then the horizontal sections of the corresponding flat bundle of
conformal blocks over $\Xn$ will have non-trivial monodromy around the
diagonals. For example, in the case of the Kac-Moody vertex algebra
$V_k(\g)$ and $X={\mathbb P}^1$, these sections are solutions of the
Knizhnik-Zamolodchikov equations, and the monodromy matrices are given
by the $R$--matrices of the quantum group $U_q(\g)$, see \cite{TK,SV}.
\end{rema}

\section{Sheaves of conformal blocks on moduli spaces}
\label{sheaves}

In the previous section we associated the space of conformal blocks to
a vertex algebra $V$, a smooth curve $X$ and a point $x$ of $X$. Now
we want to understand the behavior of these spaces as we move $x$
along $X$ and vary the complex structure on $X$. More precisely, we
wish to organize the spaces of conformal blocks into a sheaf on the
moduli space ${\mf M}_{g,1}$ of smooth pointed curves of genus $g$
(here and below by moduli space we mean the corresponding moduli stack
in the smooth topology; $\D$--modules on algebraic stacks are defined
in \cite{BB,BD}). Actually, for technical reasons we prefer to work
with the spaces of coinvariants (also called covacua), which are dual
to the spaces of conformal blocks. One can construct in a
straightforward way a quasicoherent sheaf on ${\mf M}_{g,1}$, whose
fiber at $(X,x)$ is the space of coinvariants attached to $(X,x)$. But
in fact this sheaf also carries a structure of (twisted) ${\mc
D}$--module on ${\mf M}_{g,1}$, i.e., we can canonically identify the
(projectivizations of) the spaces of coinvariants attached to
infinitesimally nearby points of ${\mf M}_{g,1}$ (actually, this
$\D$--module always descends to ${\mf M}_{g}$). The key fact used in
the proof of this statement is the ``Virasoro uniformization'' of
${\mf M}_{g,1}$. Namely, the Lie algebra $\on{Der} \K$ acts
transitively on the moduli space of triples $(X,x,z)$, where $(X,x)$
are as above and $z$ is a formal coordinate at $x$; this action is
obtained by ``gluing'' $X$ from $D_x$ and $X\bs x$ (see \thmref{uni}).

In the case when the vertex algebra $V$ is rational, the corresponding
twisted ${\mc D}$--module on ${\mf M}_g$ is believed to be coherent,
i.e., isomorphic to the sheaf of sections of a vector bundle of finite
rank equipped with a projectively flat connection (this picture was
first suggested by Friedan and Shenker \cite{FrS}). This is known to
be true in the case of the Kac-Moody vertex algebra $L_k(\g)$
\cite{TUY} (see also \cite{So}) and the minimal models of the Virasoro
algebra \cite{BFM}. Moreover, in those cases the ${\mc D}$--module on
${\mf M}_g$ can be extended to a ${\mc D}$--module with regular
singularities on the Deligne-Mumford compactification $\ol{\mf
M}_{g}$, and the dimensions of the fibers at the boundary are equal to
those in ${\mf M}_{g}$. This allows one to compute these dimensions by
``Verlinde formula'' from the dimensions of the spaces of coinvariants
attached to ${\mathbb P}^1$ with three points. The latter numbers,
called the fusion rules, can in turn be found from the matrix of the
modular transformation $\tau \mapsto -1/\tau$ acting on the space of
characters (see \thmref{zhu}). The same pattern is believed to hold
for other rational vertex algebras, but as far as I know this has not
been proved in general.

The construction of the ${\mc D}$--module structure on the sheaf of
coinvariants is a special case of the general formalism of
localization of modules over Harish-Chandra pairs, due to
Beilinson-Bernstein \cite{BB}, which we now briefly recall.

\subsection{Generalities on localization.}    \label{gener}

A {\em Harish--Chandra pair} is a pair $({\mf g},K)$ where $\mf g$ is
a Lie algebra, $K$ is a Lie group, such that ${\mf k}=\on{Lie} K$ is
embedded into ${\mf g}$, and an action $Ad$ of $K$ on $\mf g$
compatible with the adjoint action of $K$ on ${\mf k}\subset{\mf g}$
and the action of $\mf k$ on $\mf g$.

A $({\mf g},K)$--{\em action} on a scheme $Z$ is the data of an
action of $\mf g$ on $Z$ (that is, a homomorphism $\rho: \mf g \to
\Theta_Z$), together with an action of $K$ on $Z$, such that (1) the
differential of the $K$--action is the restriction of the $\mf
g$--action to $\mf k$, and (2) $\rho(Ad_k(a))=k\rho(a)k\inv.$ A $({\mf
g},K)$--action is called {\em transitive} if the map ${\mf g}\otimes
{\cali O}_Z \to \Theta_Z$ is surjective, and {\em simply transitive}
if this map is an isomorphism.

A $({\mf g},K)$--{\em structure} on a scheme $S$ is a principal
$K$--bundle $\pi:Z\to S$ together with a simply transitive $({\mf
g},K)$--action on $Z$ which extends the fiberwise action of $K$.

An example of $({\mf g},K)$--structure is the homogeneous space
$S=G/K$, where $G$ is a finite--dimensional Lie group, ${\mf
g}=\on{Lie} G$, $K$ is a Lie subgroup of $G$, and we take the obvious
right action of $({\mf g},K)$ on ${\wh X}=G$. In the
finite--dimensional setting, any space with a $(\g,K)$--structure is
locally of this form. However, in the infinite-dimensional setting
this is no longer so. An example is the
$(\on{Der}\Oo,\AutO)$--structure $Z=\wh{X} \to X$ on a smooth curve
$X$, which we used above. This example can be generalized to the case
when $S$ is an arbitrary smooth scheme. Define $\wh{S}$ to be the
scheme of pairs $(x,\vec{t}_x)$, where $x\in S$ and $\vec{t}_x$ is a
system of formal coordinates at $x$. Set $\g=
\on{Der}\C[[z_1,\dots,z_n]], K= \on{Aut}\C[[z_1,\dots,z_n]]$. Then
$\wh{S}$ is a $(\g,K)$--structure on $S$, first considered by
Gelfand, Kazhdan and Fuchs \cite{GKF}.

Let $V$ be a $(\g,K)$--module, i.e., a vector space together with
actions of $\g$ and $K$ satisfying the obvious compatibility
condition.  Then we define a flat vector bundle $\V$ on any variety
$S$ with a $(\g,K)$--structure $Z$. Namely, as a vector bundle, $\V =
Z \underset{K}\times V$. Since by assumption the $\g$--action on $Z$
is simply transitive, it gives rise to a flat connection on the
trivial vector bundle $Z \times V$ over $Z$, which descends to $\V$,
because the actions of $K$ and $\g$ are compatible. In the special
case of the $(\on{Der}\Oo,\AutO)$--structure $\wh{X}$ over a smooth
curve $X$ we obtain the flat bundle $\V$ on $X$ from \secref{global
construction}.

Now consider the case when the $\g$--action is transitive but not
simply transitive (there are stabilizers). Then we can still construct
a vector bundle $\V$ on $S$ equipped with a $K$--equivariant action of
the Lie algebroid $\g_Z = \g \otimes \Oo_Z$, but we cannot obtain an
action of $\Theta_Z$ on $\V$ because the map $a: \g_Z \arr \Theta_{Z}$
is no longer an isomorphism. Nevertheless, $\Theta_Z$ will act on any
sheaf, on which $\on{Ker} a$ acts by $0$.
In particular, the sheaf $V\otimes\Oo_{Z}/\on{Ker} a \cdot
(V\otimes\Oo_{Z}) \simeq \D_{Z}\otimes_{U\g} V$ gives rise to a
$K$--equivariant $\D$--module on $Z$, and hence to a $\D$--module on
$S$, denoted $\Delta(V)$ and called the {\em localization} of $V$ on
$S$.

More generally, suppose that $V$ is a module over a Lie algebra $\lf$,
which contains $\g$ as a Lie subalgebra and carries a compatible
$K$--action. Suppose also that we are given a $K$--equivariant Lie
algebra subsheaf $\wt\lf$ of the constant sheaf of Lie algebras $\lf
\otimes \Oo_Z$, which is preserved by the natural action of the Lie
algebroid $\g_Z$. Then if $\wt{\lf}$ contains $\on{Ker} a$, the sheaf
$V\otimes\Oo_{Z}/\wt{\lf} \cdot (V\otimes\Oo_{Z})$ is a
$K$--equivariant ${\mc D}$--module on $Z$, which descends to a ${\mc
D}$--module $\wt{\Delta}(V)$ on $S$.

Let us describe the fibers of $\wt{\Delta}(V)$ (considered as an
$\Oo_S$--module). For $s \in S$, let $Z_s$ be the fiber of $Z$ at $s$,
and $\lf_s = Z_s \underset{K}{\times} \lf$, $\V_s = Z_s
\underset{K}{\times} V$. The fibers of $\wt{\lf}$ at the points of
$Z_s \subset M$ give rise to a well-defined Lie subalgebra
$\wt{\lf}_s$ of $\lf_s$. Then the fiber of $\wt{\Delta}(V)$ at $s \in
S$ is canonically isomorphic to the space of coinvariants
$\V_s/\wt{\lf}_s \cdot \V_s$.

\subsection{Localization on the moduli space of curves.}

We apply the above formalism in the case when $S$ is the moduli space
${\mf M}_{g,1}$ of smooth pointed curves of genus $g>1$, and $Z$ is
the moduli space $\wh{\mf M}_{g,1}$ of triples $(X,x,z)$, where $(X,x)
\in {\mf M}_{g,1}$ and $z$ is a formal coordinate at $x$. Clearly,
$\wh{\mf M}_{g,1}$ is an $\AutO$--bundle over ${\mf M}_{g,1}$.

\begin{theo}[cf. \cite{ADKP,BS,Ko,TUY}]    \label{uni}
$\wh{\mf M}_{g,1}$ carries a transitive action of $\on{Der} \K$
compatible with the $\AutO$--action along the fibers.
\end{theo}

The action of
the corresponding ind--group $\on{Aut} \K$ is defined by
the ``gluing'' construction. If $(X,x,z)$ is an $R$--point of $\wh{\mf
M}_{g,1}$, and $\rho\in\on{Aut}\K(R)$, we construct a new $R$--point
$(X_\rho,x_\rho,z_\rho)$ of $\wh{\mf M}_{g,1}$ by ``gluing'' the
formal neighborhood of $x$ in $X$ with $X\bs x$ with a ``twist'' by
$\rho$.




Now we are in the situation of \secref{gener}, with $\g = \on{Der} \K$
and $K = \AutO$. Denote by ${\mc A}$ the Lie algebroid $\on{Der} \K
\wh{\otimes} \Oo_{\wh{\mf M}_{g,1}}$ on $\wh{\mf M}_{g,1}$. By
construction, the kernel of the corresponding homomorphism $a: {\mc A}
\to \Theta$ is the subsheaf ${\mc A}_{\out}$ of ${\mc A}$, whose fiber
at $(X,x) \in {\mf M}_{g,1}$ is $\on{Vect}(X\bs x)$. Applying the
construction of \secref{gener} to $V = \Vir_0$ (the Virasoro vertex
algebra with $c=0$), we obtain a ${\mc D}$--module $\Delta(\Vir_0)$ on
${\mf M}_{g,1}$, whose fibers are the spaces of coinvariants
$\Vir_0/\on{Vect}(X\bs x) \cdot \Vir_0$.

More generally, let $V$ be an arbitrary conformal vertex algebra with
central charge $0$. Then it is a module over the Lie algebra $\lf =
U(V)$ from \secref{coinv}. Let $\wt{\lf} = {\mc U}(V)_{\out}$ be the
subsheaf of $U(V) \wh{\otimes} \Oo_{\wh{\mf M}_{g,1}}$, whose fiber at
$(X,x,z) \in {\mf M}_{g,1}$ equals $U_{X\bs x}(V)$ (see
\secref{coinv}). Note that $U(V)$ contains $\on{Der} \K$, and ${\mc
U}(V)_{\out}$ contains ${\mc A}_{\out}$. Moreover, we have:

\begin{lemm}    \label{normal}
${\mc U}(V)_{\out}$ is preserved by the action of the Lie algebroid
${\mc A}$.
\end{lemm}

This is equivalent to the statement that $U_{X\bs x}(V)$ and
$U_{X_\rho\bs x_\rho}(V)$ are conjugate by $\rho$, which follows from
the fact that formula \eqref{Hformula} is valid for any $\rho \in
\on{Aut} \K$.

So we are again in the situation of \secref{gener}, and hence
we obtain a ${\mc D}$--module $\wt{\Delta}(V)$ on ${\mf M}_{g,1}$,
whose fiber at $(X,x)$ is precisely the space of coinvariants
$H(X,x,V) = \V_x/U_{X\bs x}(\V_x) \cdot \V_x$ (see \secref{coinv}),
which is what we wanted.


In the case of vertex algebras with non-zero central charge, we need
to modify the general construction of \secref{gener} as follows.
Suppose that $\g$ has a central extension $\wh\g$ which splits over
$\kk$,
and such that the extension
\begin{equation}    \label{splits}
0\to\Oo_{Z} \to\ghat \otimes \Oo_{Z} \to\g_{Z}\to 0
\end{equation}
splits over the kernel of $a:\g_Z \to \Theta_{Z}$. Then the quotient
$\g'_Z$ of $\ghat\otimes\Oo_{Z}$ by $\on{Ker}a$ is an extension of
$\Theta_{Z}$ by $\Oo_{Z}$, with a natural Lie algebroid structure. The
enveloping algebra of this Lie algebroid is a sheaf ${\mc D}'_Z$ of
twisted differential operators (TDO) on $Z$, see \cite{BB}. Further,
$\g'_Z$ descends to a Lie algebroid on $S$, and gives rise to a sheaf
of TDO ${\mc D}'_S$.

Now if $V$ is a $(\ghat,K)$--module, we obtain a $\D'_Z$--module ${\mc
D}'_Z \otimes_{U\ghat} V$. By construction, it is $K$--equivariant,
and hence descends to a $\D'_S$--module on $S$. More generally,
suppose that $\wh{\g}$ is a Lie subalgebra of $\wh{\lf}$, and we are
given a subsheaf $\wt{\lf}$ of $\wh{\lf} \otimes \Oo_Z$, containing
$\on{Ker} a$ and preserved by the action of $\wh\g \otimes
\Oo_Z$. Then the sheaf $V\otimes\Oo_{Z}/\wt{\lf} \cdot
V\otimes\Oo_{Z}$ is a $K$--equivariant ${\mc D}'_Z$--module on $Z$,
which descends to a ${\mc D}'_S$--module (still denoted by
$\wt{\Delta}(V)$). Its fibers are the coinvariants $\V_s/\wt{\lf}_s
\cdot \V_s$.

Let $V$ be a conformal vertex algebra with an arbitrary central charge
$c$. Then by the residue theorem the corresponding sequence
\eqref{splits} with $\wh\g = Vir$ splits over $\on{Ker} a = {\mc
A}_{\out}$. Hence there exists a TDO sheaf $\D_c$ on ${\mf M}_{g,1}$
and a $\D_c$--module $\wt{\Delta}(V)$ on ${\mf M}_{g,1}$ whose fiber
at $(X,x)$ is the space of coinvariants $H(X,x,V)$. Actually,
$\wt{\Delta}(V)$ descends to ${\mf M}_g$. Explicit computation shows
that the sheaf of differential operators on the determinant line
bundle over ${\mf M}_g$ corresponding to the sheaf of relative
$\la$--differentials on the universal curve is $\D_{c(\la)}$, where
$c(\la) = -12\la^2+12\la-2$ (see \cite{BS,BFM}).

It is straightforward to generalize the above construction to the case
of multiple points with arbitrary $V$--module insertions
$M_1,\ldots,M_n$. In that case we obtain a twisted $\D$--module on the
moduli space ${\mf M}'_{g,n}$ of $n$--pointed curves with non-zero
tangent vectors at the points.

\subsection{Localization on other moduli spaces.}\label{localization}

In the previous section we constructed twisted ${\mc D}$--modules on
the moduli spaces of pointed curves using the interior action of the
Harish-Chandra pair $(Vir,\AutO)$ on conformal vertex algebras and its
modules. This construction can be generalized to the case of other
moduli spaces if we consider interior actions of other Harish--Chandra
pairs.

For example, consider the Kac-Moody vertex algebra $V_k(\g)$. It
carries an action of the Harish-Chandra pair $(\ghat,G(\Oo))$, where
$G$ is the simply-connected group with Lie algebra $\g$. The
corresponding moduli space is the moduli space ${\mf M}_G(X)$ of
$G$--bundles on a curve $X$. For $x \in X$, let $\wh{\mf M}_G(X)$ be
the moduli spaces of pairs $({\mc P},s)$, where ${\mc P}$ is a
$G$--bundle on $X$, and $s$ is a trivialization of ${\mc
P}|_{D_x}$. Then $\wh{\mf M}_G(X)$ is a principal $G(\Oo_x)$--bundle
over ${\mf M}_G(X)$, equipped with a compatible
$\g(\K_x)$--action. The latter is given by the gluing construction
similar to that explained in the proof of \thmref{uni}.

Now we can apply the construction of \secref{gener} to the
$(\ghat_x,G(\Oo_x))$--module $\V_k(\g)_x$, and obtain a twisted
$\D$--module on ${\mf M}_G(X)$. Its fiber at ${\mc P} \in {\mf
M}_G(X)$ is the space of coinvariants $\V_k(\g)_x/\g_{x,\out}^{\mc P}
\cdot \V_k(\g)_x$, where $\g_{x,\out}^{\mc P}$ is the Lie algebra of
sections of ${\mc P} \underset{G}\times \g$ over $X\bs x$.

More generally, let $V$ be a vertex algebra equipped with a
$\ghat$--structure of level $k\neq -h^\vee$ (see \secref{gentw}). Then
$V$ carries an action of the Harish-Chandra pair
$(\ghat_x,G(\Oo_x))$. In the same way as in the case of the Virasoro
algebra, we attach to $V$ a twisted ${\mc D}$--module on ${\mf
M}_G(X)$, whose fiber at ${\mc P}$ is the twisted space of
coinvariants $H^{\mc P}(X,x,V)$, dual to the space $C^{\mc P}(X,x,V)$
defined in \remref{twconf}.

For the Heisenberg vertex algebra $\pi$, the corresponding moduli
space is the Jacobian $J(X)$ of $X$. The bundle $\wh{J}(X)$ over
$J(X)$, which parameterizes line bundles together with trivializations
on $D_x$, carries a transitive action of the abelian Lie algebra
$\K$. Applying the above construction, we assign to any conformal
vertex algebra equipped with an embedding $\pi \to V$ a twisted
$\D$--module on $J(X)$, whose fiber at $\Ll \in J(X)$ is the twisted
space of coinvariants $H^{\Ll}(X,x,V)$. The corresponding TDO sheaf is
the sheaf of differential operators acting on the theta line bundle on
$J(X)$.

By considering the action of the semi-direct product $Vir \ltimes
\ghat$ on $V$ we obtain a twisted $\D$-module on the moduli space of
curves and $G$--bundles on them. We can further generalize the
construction by inserting modules at points of the curve, etc.

\subsection{Local and global structure of moduli spaces.}
\label{globalstr}

In the simplest cases the twisted $\D$--modules obtained by the above
localization construction are the twisted sheaves $\D$ themselves
(considered as left modules over themselves). For example, the sheaf
$\Delta(V_k(\g))$ on ${\mf M}_G(X)$ is the sheaf $\D_k$ of
differential operators acting on the $k$th power of the determinant
line bundle ${\mc L}_G$ on ${\mf M}_G(X)$ (the ample generator of the
Picard group of ${\mf M}_G(X)$). The dual space to the stalk of $\D_k$
at ${\mc P} \in {\mf M}_G(X)$ is canonically identified with the space
of sections of ${\mc L}^{\otimes k}_G$ on the formal neighborhood of
${\mc P}$ in ${\mf M}_G(X)$. But according to our construction, this
space is also isomorphic to the space of conformal blocks $C^{\mc
P}(X,x,V_k(\g))$. In particular, the coordinate ring of the formal
deformation space of a $G$--bundle ${\mc P}$ on a curve $X$ is
isomorphic to $C^{\mc P}(X,x,V_0(\g))$. Thus, using the description of
conformal blocks in terms of correlation functions (see
\secref{functreal}), we obtain a realization of the formal deformation
space and a line bundle on it in terms of polydifferentials on powers
of the curve $X$ satisfying bootstrap conditions on the diagonals (see
\cite{BG,BD,Gi}). On the other hand, if we replace $X$ by $D_x$, the
corresponding space of polydifferentials becomes $V_k(\g)^*$ (more
precisely, its twist by the torsor of formal coordinates at
$x$). Therefore the vertex algebra $V_k(\g)$ may be viewed as the
local object responsible for deformations of $G$--bundles on curves.

Similarly, the Virasoro vertex algebra $\Vir_c$ and the Heisenberg
vertex algebra $\pi$ may be viewed as the local objects responsible
for deformations of curves and line bundles on curves, respectively.

Optimistically, one may hope that any one--parameter family of ``Verma
module type'' vertex algebras (such as $\Vir_c$ or $V_k(\g)$) is
related to a moduli space of curves with some additional
structures. It would be interesting to identify the deformation
problems related to the most interesting examples of such vertex
algebras (for instance, the $\W$--algebras of \secref{Walg}), and to
construct the corresponding global moduli spaces. One can attach a
``formal moduli space'' to any vertex algebra $V$ by taking the double
quotient of the ind--group, whose Lie algebra is $U(V)$, by its ``in''
and ``out'' subgroups. But this moduli space is very big, even in the
familiar examples of Virasoro and Kac-Moody algebras. The question is
to construct a smaller formal moduli space with a line bundle, whose
space of sections is isomorphic to the space of conformal blocks
$C(X,x,V)$.

While ``Verma module type'' vertex algebras are related to the local
structure of moduli spaces, rational vertex algebras (see
\secref{rva}) can be used to describe the global structure.

For instance, the space of conformal blocks $C(X,x,L_k(\g)) \simeq
\on{Hom}_{\g_{\on{out}}}(L_k(\g),\C)$, is isomorphic to $\Gamma_k =
\Gamma({\mf M}_G(X),\Ll^{\otimes k}_G)$, see \cite{BL,Fa,KNR}. The
space $C(X,x,L_k(\g))$ satisfies a factorization property: its
dimension does not change under degenerations of the curve into curves
with nodal singularities. This property allows one to find this
dimension by computing the spaces of conformal blocks in the case of
${\mbb P}^1$ and three points with $L_k(\g)$--module insertions (fusion
rules). The resulting formula for $\dim C(X,x,L_k) = \dim \Gamma_k$ is
called the Verlinde formula, see \cite{TUY,So}.

Since $\Ll_G$ is ample, we can recover the moduli space of semi-stable
$G$--bundles on $X$ as the $\on{Proj}$ of the graded ring
$\oplus_{k\geq 0} \Gamma_{Nk} = \oplus_{k\geq 0} C(X,x,L_{Nk}(\g))$
for large $N$ (the ring of ``non-abelian theta functions''). Thus, if
we could define the product on conformal blocks in a natural way, we
would obtain a description of the moduli space. Using the correlation
function description of conformal blocks Feigin and Stoyanovsky
\cite{FS} have identified the space $C(X,x,L_k(\g))$ with the space of
sections of a line bundle on the power of $X$ satisfying certain
conditions. For example, when $\g=\sw_2$ it is the space of sections
of the line bundle $\Omega(2nx)^{\boxtimes nk}$ over $X^{nk}$, which
are symmetric and vanish on the diagonals of codimension $k$ (here $x
\in X$, and $n \geq g$). In these terms the multiplication $\Gamma_k
\times \Gamma_l \to \Gamma_{k+l}$ is just the composition of exterior
tensor product and symmetrization of these sections. Thus we obtain an
explicit description of the coordinate ring of the moduli space of
semi-stable $SL_2$--bundles. On the other hand, replacing in the above
description $X$ by $D_x$, we obtain a functional realization of
$L_k(\sw_2)^*$.

Similarly, the space of conformal blocks corresponding to the lattice
vertex superalgebra $V_{\sqrt{N}\Z}$ (see \secref{bfc}) may be
identified with the space of theta functions of order $N$ on $J(X)$,
so we obtain the standard ``functional realization'' of the Jacobian
$J(X)$.



\subsection{Critical level.}    \label{critical}

When $k=-h^\vee$ (which is called the critical level), the vertex
algebra $V_k(\g)$ is not conformal. Nevertheless, $V_{-h^\vee}(\g)$
still carries an action of $\on{Der} \Oo$ which satisfies the same
properties as for non-critical levels. Therefore we can still attach
to $V_{-h^\vee}$ a twisted $\D$--module on ${\mf M}_G(X)$ (but we
cannot vary the curve). This $\D$--module is just the sheaf $\D'$ of
differential operators acting on the square root of the canonical
bundle on ${\mf M}_G(X)$ (see \cite{BD}).

Recall from \secref{coset} that the center of a vertex algebra $V$ is
the coset vertex algebra of the pair $(V,V)$. This is a commutative
vertex subalgebra of $V$. The center of $V_k(\g)$ is simply its
subspace of $\g[[z]]$--invariants. It is easy to show that this
subspace equals $\C v_k$ if $k\neq -h^\vee$. In order to describe the
center ${\mf z}(\ghat)$ of $V_{-h^\vee}(\g)$ we need the concept of
{\em opers}, due to Beilinson and Drinfeld \cite{BD}.

Let $G_{\on{ad}}$ be the adjoint group of $\g$, and $B_{\on{ad}}$ be
its Borel subgroup. By definition, a $\g$--oper on a smooth curve $X$
is a $G_{\on{ad}}$--bundle on $X$ with a (flat) connection $\nabla$
and a reduction to $B_{\on{ad}}$, satisfying a certain transversality
condition \cite{BD}. For example, an $\sw_2$--oper is the same as a
projective connection. The set of $\g$--opers on $X$ is the set of
points of an affine space, denoted $\on{Op}_{\g}(X)$, which is a
torsor over $\oplus_{i=1}^\ell \Gamma(X,\Omega^{d_i+1})$, where
$d_i$'s are the exponents of $\g$. Denote by $A_{\g}(X)$ the ring of
functions on $\on{Op}_{\g}(X)$.

The above definition can also be applied when $X$ is replaced by the
disc $D = \on{Spec} \C[[z]]$ (see \cite{DS}). The corresponding ring
$A_{\g}(D)$ has a natural $(\on{Der} \Oo,\AutO)$--action, and is
therefore a commutative vertex algebra (see \secref{commex}). It is
isomorphic to a limit of the $\W$--vertex algebra $\W_k(\g)$ as $k \to
\infty$ (see \secref{Walg}), and because of that it is called the
classical $\W$--algebra corresponding to $\g$. For example,
$A_{\sw_2}(D)$ is a limit of $\on{Vir}_c$ as $c \to \infty$. The
following result was conjectured by Drinfeld and proved in
\cite{FF:gd}.

\begin{theo}    \label{center}
The center ${\mf z}(\ghat)$ of $V_{-h^\vee}(\g)$ is isomorphic, as a
commutative vertex algebra with $\on{Der} \Oo$--action, to
$A_{^L\g}(D)$, where $^L\g$ is the Langlands dual Lie algebra to $\g$.
\end{theo}

In fact, more is true: both of these commutative vertex algebras carry
what can be called Poisson vertex algebra (or coisson algebra, in the
terminology of \cite{BD:ch}) structures, which are preserved by this
isomorphism.

The localization of ${\mf z}(\ghat)$ on ${\mf M}_G(X)$ is the
$\Oo$--module $H(X,x,{\mf z}(\ghat)) \otimes \Oo_{{\mf M}_G(X)}$
(since $\ghat$ does not preserve ${\mf z}(\ghat) \subset
V_{-h^\vee}(\g)$, this sheaf does not carry a $\D$--module
structure). By functoriality, the embedding ${\mf z}(\ghat) \to
V_{-h^\vee}(\g)$ of vertex algebras gives rise to a homomorphism of
$\Oo$--modules $H(X,x,{\mf z}(\ghat)) \otimes \Oo_{{\mf M}_G(X)} \to
\D'$. \thmref{center} implies that $H(X,x,{\mf z}(\ghat)) \simeq
H(X,x,A_{^L\g}(D))$. As shown in \cite{BD:ch}, the space of
coinvariants of a commutative vertex algebra $A$ is canonically
isomorphic to the ring of functions on $\Gamma_\nabla(X,\on{Spec} {\mc
A})$. In our case we obtain: $\Gamma_\nabla(X,\on{Spec} {\mc
A}_{^L\g}) \simeq \on{Op}_{^L\g}(X)$, and so $H(X,x,A_{^L\g}(D))
\simeq A_{^L\g}(X)$. Passing to global sections we obtain:

\begin{coro}
There is a homomorphism of algebras $A_{^L\g}(X) \to \Gamma({\mf
M}_G(X),\D')$.
\end{coro}

Beilinson and Drinfeld have shown in \cite{BD} that this map is an
embedding, and if $G$ is simply-connected, it is an isomorphism. In
that case, each $^L G$--oper $\rho$ gives us a point of the spectrum
of the commutative algebra $\Gamma({\mf M}_G(X),\D')$, and hence a
homomorphism $\wt{\rho}: \Gamma({\mf M}_G(X),\D') \to \C$. The
$\D'$--module $\D'/(\D' \cdot m_{\wt{\rho}})$ on ${\mf M}_G(X)$, where
$m_{\wt{\rho}}$ is the kernel of $\wt{\rho}$, is holonomic. It is
shown in \cite{BD} that the corresponding untwisted $\D$--module is a
Hecke eigensheaf attached to $\rho$ by the geometric Langlands
correspondence.

\subsection{Chiral deRham complex.}    \label{cdo}

Another application of the general localization pattern of
\secref{gener} is the recent construction, due to Malikov, Schechtman
and Vaintrob \cite{MSV}, of a sheaf of vertex superalgebras on any
smooth scheme $S$ called the chiral deRham complex. Under certain
conditions, this sheaf has a purely even counterpart, the sheaf of
chiral differential operators. Here we give a brief review of the
construction of these sheaves.

Let $\Gamma_N$ be the Weyl algebra associated to the symplectic vector
space $\C((t))^N \oplus (\C((t)) dt)^N$. It has topological generators
$a_{i,n}, a^*_{i,n}, i=1,\ldots,N; n \in \Z$, with relations
$$[a_{i,n},a^*_{j,m}]=\delta_{i,j} \delta_{n,-m}, \quad \quad
[a_{i,n},a_{j,m}] = [a^*_{i,n},a^*_{j,m}] = 0.
$$
Denote by $H_N$ the Fock representation of $\Gamma_N$ generated by the
vector $\vac$, satisfying $a_{i,n} \vac = 0, n\geq 0; a^*_{i,m} \vac =
0, m>0$. This is a vertex algebra, such that $$Y(a_{i,-1}\vac,z) =
a(z) = \sum_{n\in\Z} a_{i,n} z^{-n-1}, \quad \quad Y(a^*_{i,0}\vac,z)
= a^*(z) = \sum_{n\in\Z} a^*_{i,n} z^{-n}.$$ Let $\bigwedge^\bullet_N$
be the fermionic analogue of $H_N$ defined in the same way as in
\secref{Walg}. This is the Fock module over the Clifford algebra with
generators $\psi_{i,n}, \psi^*_{i,n}, i=1,\ldots,N, n\in\Z$, equipped
with an additional $\Z$--gradation. It is shown in \cite{MSV} that the
vertex algebra structure on $H_N$ may be extended to $\wh{H}_N = H_N
\us{\C[a^*_{i,0}]}\otimes \C[[a^*_{i,0}]]$. Let $\wh{\Omega}^\bullet_N
= \wh{H}_N \otimes \bigwedge^\bullet_N$. It carries a differential $d
= \on{Res} \sum_i a_i(z) \psi^*_i(z)$, which makes it into a complex.

Denote by $W_N$ (resp., $\Oo_N$, $\Omega^1_N$) the topological Lie
algebra of vector fields (resp., ring of functions, module of
differentials) on $\on{Spec} \C[[t_1,\ldots,t_N]]$. Define a map $W_N
\to U(\wh{\Omega}_N^\bullet)$ (the completion of the Lie algebra of
Fourier coefficients of fields from $\wh{\Omega}_N^\bullet$) sending
$$
f(t_i) \pa_{t_j} \mapsto \on{Res} \left( :f(a^*_i(z)) a_j(z): +
\sum_{k=1}^N :(\pa_{t_k} f)(a_i^*(z)) \psi^*_k(z) \psi_j(z): \right).
$$
Explicit computation shows that it is a homomorphism of Lie algebras
(this is an example of ``supersymmetric cancellation of anomalies'')
\cite{FF:si,MSV}. We obtain a $W_N$--action on
$\wh{\Omega}^\bullet_N$, which can be exponentiated to an action of the
Harish-Chandra pair $(W_N,\AutO_N)$. Furthermore, because the
$(W_N,\AutO_N)$--action comes from the residues of vertex operators,
it preserves the vertex algebra structure on $V$ (so the group
$\AutO_N$ is an example of the group of internal symmetries from
\secref{gentw}). Recall from \secref{gener} that any smooth scheme $S$
of dimension $N$ carries a canonical $(W_N,\AutO_N)$--structure
$\wh{S}$. Applying the construction of \secref{gener}, we attach to
$\wh{\Omega}^\bullet_N$ a $\D$--module on $S$. The sheaf of horizontal
sections of this $\D$--module is a sheaf of vertex superalgebras (with
a structure of complex) on $S$. This is the {\em chiral deRham
complex} of $S$, introduced in \cite{MSV}. For its applications to
mirror symmetry and elliptic cohomology, see \cite{BoL,Bor1}.

Now consider a purely bosonic analogue of this construction (i.e.,
replace $\wh{\Omega}_N$ with $\wh{H}_N$). The Lie algebra
$U(\wh{H}_N)$ of Fourier coefficients of the fields from $\wh{H}_N$
has a filtration $A^{\leq i}_N, i\geq 0$, by the order of
$a_{i,n}$'s. Note that $A^{\leq 1}_N$ is a Lie subalgebra of
$U(\wh{H}_N)$ and $A^0_N$ is its (abelian) ideal. Denote by $\T_N$ the
Lie algebra $A^{\leq 1}_N/A^{0}_N$. Intuitively, $U(\wh{H}_N)$ is the
algebra of differential operators on the topological vector space
$\C((t))^N$ with the standard filtration, $A^{(0)}$ (resp., $\T_N$) is
the space of functions (resp., the Lie algebra of vector fields) on
it, and $\wh{H}_N$ is the module of ``delta--functions supported on
$\C[[t]]^N$''.

In contrast to the finite-dimensional case, the extension $$0 \to
A^{0}_N \to A_N^{\leq 1} \to \T_N \to 0$$ does not split \cite{FF:si}.
Because of that, the naive map $W_N \to U(\wh{H}_N)$ sending $f(t_i)
\pa_{t_j}$ to $\on{Res} :f(a^*_i(z)) a_j(z):$ is not a Lie algebra
homomorphism. Instead, we obtain an extension
\begin{equation}    \label{WN}
0 \to \Omega^1_N/d\Oo_N \to \wt{W}_N \to W_N \to 0
\end{equation}
of $W_N$ by its module $\Omega^1_N/d\Oo_N$, which is embedded into
$A^{\leq 1}_N$ by the formula $g(t_i) dt_j \mapsto \on{Res}
g(a^*_i(z)) \pa_z a^*_j(z)$. If the sequence \eqref{WN} were split, we
would obtain a $(W_N,\AutO_N)$--action on $\wh{H}_N$ and associate to
$\wh{H}_N$ a $\D$--module on $S$ in the same way as above.

But in reality the extension \eqref{WN} does not split for $N>1$ (note
that $\Omega^1_1/d\Oo_1 = 0$). Hence in general we only have an action
on $\wh{H}_N$ of the Harish-Chandra pair $(\wt{W}_N,\wt{\on{Aut}}
\Oo_N)$, where $\wt{\on{Aut}} \Oo_N$ is an extension of $\AutO_N$ by
the additive group $\Omega^1_N/d\Oo_N$. This action again preserves
the vertex algebra structure on $V$. To apply the construction of
\secref{gener}, we need a lifting of the $(W_N,\AutO_N)$--structure
$\wh{S}$ on $S$ to a $(\wt{W}_N,\wt{\on{Aut}} \Oo_N)$--structure. Such
liftings form a {\em gerbe} on $S$ (in the analytic topology) with the
lien $\Omega_S^1/d\Oo_S \simeq \Omega^2_{S,\on{cl}}$, the sheaf of
closed holomorphic two-forms on $S$. The equivalence class of this
gerbe in $H^2(S,\Omega^2_{S,\on{cl}})$ is computed in \cite{GMS}.
If this class equals $0$, then we can lift $\wh{S}$ to an
$(\wt{W}_N,\wt{\on{Aut}} \Oo_N)$--structure $\wt{S}$ on $S$, and the
set of isomorphism classes of such liftings is an
$H^1(S,\Omega^2_{S,\on{cl}})$--torsor. Applying the construction of
\secref{gener}, we can then attach to $\wh{H}_N$ a $\D$--module on
$S$. The sheaf of horizontal sections of this $\D$--module is a sheaf
of vertex algebras on $S$. This is the sheaf of {\em chiral
differential operators} on $S$ (corresponding to $\wt{S}$) introduced
in \cite{MSV,GMS}. As shown in \cite{GMS}, the construction becomes
more subtle in Zariski topology; in this case one naturally obtains a
gerbe with the lien $\Omega_S^2 \to \Omega^3_{S,\on{cl}}$.

In the case when $S$ is the flag manifold $G/B$ of a simple Lie
algebra $\g$, one can construct the sheaf of chiral differential
operators in a more direct way. The action of $\g$ on the flag
manifold $G/B$ gives rise to a homomorphism $\g((t)) \to \T_N$, where
$N = \dim G/B$. It is shown in \cite{W,FF:usp,FF:si} that it can be
lifted to a homomorphism $\ghat \to A^{\leq 1}_N$ of level $-h^\vee$
(Wakimoto realization). Moreover, it gives rise to a homomorphism of
vertex algebras $V_{-h^\vee}(\g) \to H_N$. Hence we obtain an
embedding of the constant subalgebra $\g$ of $\ghat$ into $\wt{W}_N$,
and a $(\g,B)$--action on $\wh{H}_N$. On the other hand, $G/B$ has a
natural $(\g,B)$--structure, namely $G$. Applying the localization
construction to the Harish-Chandra pair $(\g,B)$ instead of
$(\wt{W}_N,\wt{\on{Aut}} \Oo_N)$, we obtain the sheaf of chiral
differential operators on $G/B$ -- note that it is uniquely defined in
this case (see \cite{MSV}).


\section{Chiral algebras}    \label{Chiral Algebras}

In \secref{global construction}, we gave a coordinate independent
description of the vertex operation $Y(\cdot,z)$. The formalism of
chiral algebras invented by Beilinson and Drinfeld \cite{BD:ch} (see
\cite{Ga} for a review) is based on a coordinate--free realization of
the operator product expansion (see \secref{associativ}). In the
definition of $\Y_x$ we acted by $Y(A,z)$ on $B$, placed at the fixed
point $x \in X$. In order to define the OPE invariantly, we need to
let $x$ move along $X$ as well.  This is one as follows. Recall that
to each conformal vertex algebra $V$ we have assigned a vector bundle
$\V$ on any smooth curve $X$. Choose a formal coordinate $z$ at $x \in
X$, and use it to trivialize $\V|_{D_x}$. Then $\V\boxtimes
\V(\infty\Delta)$ is a sheaf on $D_x^2 = \on{Spec} \C[[z,w]]$
associated to the $\C[[z,w]]$--module $V \otimes
V[[z,w]][(z-w)^{-1}]$. Let $\Delta_! \V$ be a sheaf on $D_x^2$
associated to the $\C[[z,w]]$--module
$V[[z,w]][(z-w)^{-1}]/V[[z,w]]$. Recall from \secref{global
construction} that we have a flat connection on $\V$. Hence $\V$ is a
vector bundle with a flat connection, i.e., a ${\mc D}$--module, on
$D_x$. The sheaves $\V\boxtimes \V(\infty\Delta)$ and $\Delta_!  \V$
are not vector bundles, but they have natural structures of
$\D$--modules on $D_x^2$ (independent of the choice of $z$). Because
of that, we switch to the language of $\D$--modules. The following
result is proved in the same way as in \thmref{universal section} (see
also \cite{HL}).

\begin{theo}\label{invariant OPE}
Define a map $\Y_x^{(2)}:\V\boxtimes \V(\infty\Delta) \to \Delta_!\V$
by the formula
$$\Y_x^{(2)}(f(z,w) A \boxtimes B) = f(z,w)Y(A,z-w)\cdot B \quad
\on{mod} V[[z,w]].$$ Then $\Y_x^{(2)}$ is a homomorphism of ${\mc
D}$--modules, which is independent of the choice of the coordinate
$z$.
\end{theo}

At this point it is convenient to pass to the right $\D$--module on
$X$ corresponding to the left $\D$--module $\V$. Recall that if ${\mc
F}$ is a left $\D$ module on a smooth variety $Z$, then ${\mc F}^r :=
{\mc F} \otimes \Omega_Z$, where $\Omega_Z$ denotes the canonical
sheaf on $Z$, is a right $\D$--module on $Z$. Denote $\Delta: X
\hookrightarrow X^2, j: X^2\backslash \Delta(X) \hookrightarrow X^2$.

\begin{coro}    \label{chal}
The vertex algebra structure on $V$ gives rise to a homomorphism of
right $\D$--modules on $X^2$, $\mu: j^* j_* \V^r \boxtimes \V^r \to
\Delta_!  \V^r$ satisfying the following conditions:
\begin{enumerate}
\item[$\bullet$] {\em (skew-symmetry)} $\mu(f(x,y) A \boxtimes B) = -
\sigma_{12} \circ \mu(f(y,x) B \boxtimes A)$;
\item[$\bullet$] {\em (Jacobi identity)} $\mu(\mu(f(x,y,z) A \boxtimes
B) \boxtimes C) + \sigma_{123} \circ \mu(\mu(f(y,z,x) B \boxtimes C)
\boxtimes A)$
$$+ \sigma_{123}^{-1} \circ \mu(\mu(f(z,x,y) C \boxtimes A) \boxtimes
B) = 0.$$
\item[$\bullet$] {\em (vacuum)} we are given an embedding $\Omega
\hookrightarrow \V^r$ compatible with the natural homomorphism $j_*
j^* \Omega \boxtimes \Omega \to \Delta_!  \Omega$.
\end{enumerate}
\end{coro}

Beilinson and Drinfeld define a {\em chiral algebra} on a curve $X$ as
a right $\D$--module ${\mc A}$ on $X$ equipped with a homomorphism $j_*
j^* {\mc A} \boxtimes {\mc A} \to \Delta_! {\mc A}$ satisfying the
conditions of \corref{chal} (see \cite{BD:ch,Ga}).

The axioms of chiral algebra readily imply that for any chiral algebra
${\mc A}$, its deRham cohomology $H_{\on{dR}}^\bullet(\Sigma,{\mc A})$
is a (graded) Lie algebra for any open $\Sigma \subset X$.
In particular, for an affine open $\Sigma \subset X$ we obtain the
result mentioned in \secref{coinv} that $H_{\on{dR}}^0(\Sigma,\V^r)$
is a Lie algebra.

Define the right $\D$--module $\V_2^r$ on $X^2$ as the kernel of the
homomorphism $\mu$, and let $\V_2 = \V_2^r \otimes \Omega_{X^2}^{-1}$
be the corresponding left $\D$--module. There is a canonical
isomorphism $H(X,x,V) \simeq H_{\on{dR}}^2(X^2,\V_2^r)$ (see
\cite{Ga}, Prop. 5.1). It is dual to the isomorphism $C(X,x,V) \simeq
\Gamma_\nabla(X^2,\V_2^*)$ mentioned in \secref{functreal}.

Beilinson and Drinfeld give a beautiful description of chiral algebras
as {\em factorization algebras}.
By definition, a factorization algebra is a collection of
quasicoherent $\Oo$--modules ${\mc F}_n$ on $X^n, n\geq 1$, satisfying
a factorization condition, which essentially means that the fiber of
${\mc F}_n$ at $(x_1,\ldots,x_n)$ is isomorphic to $\otimes_{s \in S}
({\mc F}_1)_{s}$, where $S=\{x_1,\ldots,x_n\}$. For example, $j^* {\mc
F}_2 = j^* ({\mc F}_1 \boxtimes {\mc F}_1)$ and $\Delta^* {\mc F}_2 =
{\mc F}_1$. Intuitively, such a collection may be viewed as an
$\Oo$--module on the {\em Ran space} ${\mc R}(X)$ of all finite
non-empty subsets of $X$.
In addition, it is required that ${\mc F}_1$ has a global section
(unit) satisfying natural properties. It is proved in \cite{BD:ch}
that under these conditions each ${\mc F}_n$ is automatically a left
$\D$--module. Moreover, the right $\D_X$--module ${\mc F}^r_1 = {\mc
F}_1 \otimes \Omega_X$ acquires a canonical structure of chiral
algebra, with $\mu$ given by the composition $j^* j_* {\mc F}^r_1
\boxtimes {\mc F}^r_1 = j^* j_* {\mc F}^r_2 \to \Delta_! \Delta^! {\mc
F}^r_2 = \Delta_! {\mc F}^r$. In fact, there is an equivalence between
the categories of factorization algebras and chiral algebras on $X$
\cite{BD}.

For a chiral algebra $\V^r$, the corresponding sheaf ${\mc F}_n$ on
$X^n$ is the sheaf $\V_n$ mentioned in \secref{functreal} (its dual is
the sheaf of chiral correlation functions). Using these sheaves,
Beilinson and Drinfeld define ``chiral homology''
$H^{\on{ch}}_i(X,\V^r), i\geq 0$, of a chiral algebra
$\V^r$. Intuitively, this is (up to a change of cohomological
dimension) the deRham cohomology of the factorization algebra
corresponding to $\V$, considered as a $\D$--module on the Ran space
${\mc R}(X)$. The $0$th chiral homology of $\V^r$ is isomorphic to
$H_{\on{dR}}^n(X^n,\V_n \otimes \Omega_{X^n})$ for all $n>1$ and hence
to the space of coinvariants $H(X,x,V)$. The full chiral homology
functor may therefore be viewed as a derived coinvariants functor.

As an example, we sketch the Beilinson-Drinfeld construction of the
factorization algebra corresponding to the Kac-Moody vertex algebra
$V_0(\g)$ \cite{BD} (see also \cite{Ga}). Let $Gr_n$ be the ind-scheme
over $X^n$ whose fiber $Gr_{x_1,\ldots,x_n}$ over $(x_1,\ldots,x_n)
\in X$ is the moduli space of pairs $({\mc P},t)$, where ${\mc P}$ is
a $G$--bundle on $X$, and $t$ is its trivialization on $X\bs \{
x_1,\ldots,x_n \}$ (e.g., the fiber of $Gr_1$ over any $x \in X$ is
isomorphic to the affine Grassmannian $G(\K)/G(\Oo)$). Let $e$ be the
section of $Gr_n$ corresponding to the trivial $G$--bundle, and ${\mc
A}_{x_1\ldots,x_n}$ be the space of delta--functions on
$Gr_{x_1,\ldots,x_n}$ supported at $e$. These are fibers of a left
$\D$--module ${\mc A}_n$ on $X^n$. Claim: $\{ {\mc A}_n \}$ form a
factorization algebra. The factorization property for $\{ {\mc A}_n
\}$ follows from the factorization property of the ind-schemes $Gr_n$:
namely, $Gr_{x_1,\ldots,x_n} = \prod_{s \in S} Gr_s$, where
$S=\{x_1,\ldots,x_n\}$. The chiral algebra corresponding to the
factorization algebra $\{ {\mc A}_n \}$ is $\V_0(\g)^r$. To obtain the
factorization algebra corresponding to $\V_k(\g)^r$ with $k\neq 0$,
one needs to make a twist by a line bundle on $Gr_1$.


A spectacular application of the above formalism is the
Beilinson-Drinfeld construction of the {\em chiral Hecke algebra}
associated to a simple Lie algebra $\g$ and a negative integral level
$k < -h^\vee$. The corresponding vertex algebra $S_k(\g)$ is a module
over $\ghat \times \GL$: \newline $S_k(\g) = \ds \oplus_{\la \in
{}^LP^+} M_{\la,k} \otimes V_\la^*$, where $V_\la$ is a
finite-dimensional $\GL$--module with highest weight $\la$, and
$M_{\la,k}$ is the irreducible highest weight $\ghat$--module of level
$k$ with highest weight $(k+h^\vee) \la$. Here $^LG$ is the group of
adjoint type corresponding to $^L\g$, and $^LP^+$ is the set of its
dominant weights; we identify the dual $\h^*$ of the Cartan subalgebra
of $\g$ with $\h = {}^L \h^*$ using the normalized invariant inner
product on $\g$. In particular, $S_k(\g)$ is a module over $V_k(\g) =
M_{0,k}$. Note that using results of \cite{KL} one can identify the
tensor category of $\GL$--modules with a subcategory of the category
$\Oo$ of $\ghat$--modules of level $k$, so that $V_\la$ corresponds to
$M_{\la,k}$. The analogues of $S_k(\g)$ when $\g$ is replaced by an
abelian Lie algebra with a non-degenerate inner product are lattice
vertex algebras $V_L$ from \secref{rva}.

As a chiral algebra, the chiral Hecke algebra is the right
$\D$--module ${\mc S}_k(\g)^r$ on $X$ corresponding to the twist of
$S_k(\g)$ by $\wh{X}$. The corresponding factorization algebra is
constructed analogously to the above construction of $V_k(\g)$ using
certain irreducible $\D$--modules on $Gr_1$ corresponding to $V_\la$;
see \cite{Ga} for the construction in the abelian case (it is also
possible to construct a vertex algebra structure on $S_k(\g)$
directly, using results of \cite{KL}). Moreover, for any flat
$\GL$--bundle ${\mc E}$ on $X$, the twist of ${\mc S}_k(\g)^r$ by
${\mc E}$ with respect to the $\GL$--action on $S_k(\g)$ is also a
chiral algebra on $X$. Conjecturally, the complex of sheaves on ${\mf
M}_G(X)$ obtained by localization of this chiral algebra (whose fibers
are its chiral homologies) is closely related to the automorphic
$\D$--module that should be attached to ${\mc E}$ by the geometric
Langlands correspondence.


\end{document}